\documentclass[11pt,a4wide]{article}
\usepackage[utf8]{inputenc}
\usepackage{amsmath}
\usepackage{amsfonts}
\usepackage{amssymb}
\usepackage{array,multirow,graphicx}
\usepackage{hyperref}
\usepackage{a4wide}
\usepackage{relsize}
\usepackage{tikz}
\usetikzlibrary{shapes,arrows,shadows,backgrounds}
\usepackage{rotating}
\usepackage{comment}
\usepackage[numbers,sort]{natbib}
\newcommand{\dt}{ \Delta t }
\newcommand{\R}{\mathbb{R}}
\newtheorem{lem}{Lemma}[section]

\newtheorem{remark}{Remark}[section]

\newcommand{\vy}{{\vec{y}}}
\newcommand{\vu}{{\vec{u}}}
\newcommand{\vp}{{\vec{p}}}
\newcommand{\p}{{\overline{p}}}
\newcommand{\g}{{\overline{g}}}
\newcommand{\dy}{{\Delta y}}
\newcommand{\f}{{\overline{f}}}
\newcommand{\fl}{\overline{f}}

\newcommand{\la}{\overline{\lambda}}

\begin{document}
\date{\today}
\title{Linear multistep methods for optimal control problems and applications to hyperbolic relaxation systems}
\author{
G. Albi\footnote{University of Verona, Department of Computer Science, Str. Le Grazie 15, I-37134 Verona, Italy, {\tt giacomo.albi@univr.it}} 
\and 
M. Herty\footnote{RWTH Aachen University, Templergraben 55, 52062 Aachen, Germany, {\tt herty@mathc.rwth-aachen.de}}
 \and L. Pareschi\footnote{University of Ferrara, Department of Mathematics and Computer Science, Via Machiavelli 35, I-44121 Ferrara, Italy, {\tt lorenzo.pareschi@unife.it} } } 
\maketitle
\date
\begin{abstract}
We are interested in high-order linear multistep schemes for time discretization of adjoint equations
arising within optimal control problems. First we consider optimal control problems for ordinary differential equations and show loss of accuracy for Adams-Moulton and Adams-Bashford methods, whereas BDF methods preserve high--order accuracy.
Subsequently we extend these results to semi--lagrangian discretizations of hyperbolic relaxation systems. Computational results illustrate theoretical findings.  
\end{abstract}
%%%%%%%%%%%%%%%%%%%%%%%%%%%%%%%%%%%%%%%%%%%
\\
~\\
{\bf AMS:}  35L65, 49J15, 35Q93, 65L06.

%\tableofcontents

\section{Introduction}\label{introduction}

Efficient time integration methods are important for the numerical solution of optimal control
problems governed by ordinary (ODEs) and partial differential equations (PDEs).  In order to
increase  efficiency of the solvers, by reducing the memory requirements, there is a strong interest in the development of high--order methods. 
However, direct applications of standard numerical schemes to the adjoint differential systems
of the optimal control problem may lead to order reduction problems \cite{Hager00, Sandu}.  
Besides classical applications to ODEs 
these problems gained interest recently in PDEs, in particular in the field of 
hyperbolic and kinetic equations \cite{AHJP, HertyBanda11, Banda2, HertySchleper11}. 

%such methods for applications to 
%semi--Lagrangian schemes, see e.g. \cite{GroppiRussoStracquadanio}.

In this work we focus  on high--order linear multi--step methods
for optimal control problems for ordinary differential equations as well
as for semi--Lagrangian approximations of hyperbolic and kinetic transport equations, see for example \cite{DPacta,12,15,GroppiRussoStracquadanio,16,PC,CFSW,FF98}.

Regarding the time discretization of differential equations  many results in particular on Runge--Kutta methods have been established in the past years.
Properties of Runge--Kutta  methods for use in optimal control 
have been investigated for example in
\cite{Hager00,BonnansLaurent-Varin06,Sandu2,Kaya10,LangVerwer2011,DH01,DH02,HertyPareschiSteffensen2013,schroederlang}.
In particular,  Hager \cite{Hager00} investigated order conditions  for Runge--Kutta methods applied to optimality systems. This
work has been later extended \cite{BonnansLaurent-Varin06, Kaya10,HertyPareschiSteffensen2013} and also properties of symplecticity have been studied,
see also \cite{CHV08}. Further studies of discretizations of state
and control constrained problems using Runge--Kutta methods have
been  conducted in \cite{DH01,DH02,LangVerwer2011,schroederlang} 
as well as automatic differentiation techniques 
\cite{Walther2007}. Previous results for linear multi--steps method
have been considered by Sandu in \cite{Sandu}. Therein, first--order schemes 
are discussed and stability with respect to non--uniform
temporal grids has been studied. Here, we extend the results 
to high--order adjoint discretizations as well as to problems
governed by partial differential equations.  However, we restrict ourselves to the case
of uniform temporal grids. 

In the PDE context, we will focus on hyperbolic relaxation approximations to conservation laws and relaxation type kinetic equations, \cite{CSR,Pa01}. For such problems semi--Lagrangian approximations have been proposed recently in \cite{GroppiRussoStracquadanio} in combination with Runge--Kutta and BDF methods. The main advantage of such an approach is that the relaxation operator can be treated implicitly and the 
CFL condition can be circumvented by a semi-Lagrangian formulation. We mention here also \cite{DP12} where linear multistep methods have been developed for general kinetic equations.
We consider a general linear multistep setting for semi--Lagrangian schemes to reduce the optimal control problem for the PDEs to an optimal control problem for a system of ODEs.
 
The rest of the paper is organized as follows. In Section 2 we introduce the prototype optimal control problem for ODEs and consider the case of a general linear multi-step scheme. We then study the conditions under which the time discrete optimal control problem originates the corresponding time discrete adjoint equations. We prove that Adams type methods may reduce to first order accuracy and that only BDF schemes guarantee that the discretize-then-optimize approach is equivalent to the optimize-then-discretized one.
Next, in Section 3, we consider the case of semi--Lagrangian approximation of hyperbolic relaxation systems and extend the linear multistep methods to control problems for such systems. In Section 4 with the aid of several numerical examples we show the validity of our analysis. Finally we report some concluding remarks in Section 5. 

%%%%%%%%%%%%%%%%%%%%%%%%%%%%%%%%%%%%%%%%%%%
\section{Linear multi-step methods for optimal control problems of ODEs}

We are interested in linear multi--step methods for the time integration of 
ordinary differential and partial differential equations. In order to illustrate 
the approach we consider first the following problem. 

\begin{subequations}\label{OCP}
\begin{align}
 (OCP) \qquad&         \min \;  j(y(T)) \quad \mbox{ such that }\\
                         & \dot y(t)= f(y(t),u(t)),\qquad t\in [0,T]\\
               &      y(0) =  y_0 .
  \end{align}
\end{subequations}

Related to the optimal control problem we introduce the
Hamiltonian function $H$ as \begin{equation}\label{cH}
H(y,u,p):= p ^T  f(y,u).\end{equation}
Under  appropriate conditions it is well--known  \cite{Hestenes80, Troutman1996}
that the first--order optimality conditions for (\ref{OCP})
are
\begin{subequations}\label{osOCP}
\begin{eqnarray}
 \dot y&=& ~ ~ \,H_p(y,u,p)=  f(y,u), \qquad \qquad\quad ~~ ~ y(0)\,=\, y^0\label{osOCPa}\\
 \dot p&=& - H_y(y,u,p)= - f_y(y,u)^T p, \qquad p(T)\,=\, j'(y(T))\label{osOCPb}\\
0&=& ~~H_u(y,u,p)=f_u(y,u)^T p.\label{osOCPc}
\end{eqnarray}
\end{subequations}
%%%%
we assume  $f: \mathbb{R}^n \times \mathbb{R}^m \to \mathbb{R}^n$, then,  for some integer $\kappa \geq 2$, the problem (\ref{OCP}) has a local solution $(y^*,u^*)$  in $W^{\kappa,\infty} \times W^{\kappa-1,\infty}.$ There exists an open set $\Omega \subset \mathbb{R}^n \times \mathbb{R}^m$ and $\rho>0$ such that $B_\rho( y^*(t), u^*(t) ) \subset \Omega$ for every $t\in [0,T].$ If the first $\kappa$ derivatives of $f$ and $g$ are Lipschitz continuous in $\Omega$ and the first $\kappa$ derivatives of $j$ are Lipschitz in $B_\rho(y^*(T))$,  then, there exists an associated Lagrange multiplier $p^* \in W^{\kappa,\infty}$ for which the first--order optimality conditions (\ref{osOCP}) are necessarily satisfied in $(y^*,p^*,u^*).$ Under additional coercivity assumptions on the Hamiltonian (\ref{osOCP}) those conditions are  also sufficient \cite[Section 2]{Hager00}. From now on we assume that the previous conditions are fulfilled.

\begin{figure}[tb]\center
\centering
  \pgfdeclarelayer{background}
\pgfdeclarelayer{foreground}
\pgfsetlayers{background,main,foreground}

  \begin{tikzpicture}[auto,
    %decision/.style={diamond, draw=black, thick, fill=white,
    %text width=8em, text badly centered,
    %inner sep=1pt, font=\sffamily\small},
    block_center/.style ={rectangle, draw=black, thick, fill=white,
      text width=8em, text centered,
      minimum height=4em},
    block_left/.style ={rectangle,rounded corners, draw=black, thick, fill=white!50!cyan,
      text width=13em,  text centered, minimum height=4em, inner sep=6pt},
    block_noborder/.style ={rectangle,rounded corners, draw=none, thick, fill=none,
      text width=18em, text centered, minimum height=1em},
    block_assign/.style ={rectangle,rounded corners, draw=black, thick, fill=white,
      text width=18em, text ragged, minimum height=3em, inner sep=6pt},
    block_lost/.style ={rectangle, draw=black, thick, fill=white,
      text width=16em, text ragged, minimum height=3em,fill=cyan, inner sep=6pt},
      line/.style ={draw, thick,-triangle 45,fill=cyan, shorten >=6pt},
      line2/.style ={draw, thick, triangle 45-triangle 45,fill=cyan,  shorten >=0.2pt},
      description/.style={fill=white,inner sep=2pt},
        line3/.style ={draw,dashed, thick, triangle 45-triangle 45,fill=cyan,  shorten >=0.2pt},
      description/.style={fill=white,inner sep=2pt}]
    % outlining the flowchart using the PGF/TikZ matrix funtion
    \matrix [column sep=40mm,row sep=10mm] {
      % enrollment - row 1
      \node [block_left] (referred) {Continuous Optimal\\ Control Problem};
      & \node [block_left] (excluded1) {Continuous Adjoint\\ Equations}; \\
      % enrollment - row 2
     % \begin{pgfonlayer}{background}
      \node [block_left] (assessment) {BDF discretized\\ Control Problem}; 
      & \node [block_left] (excluded2) {BDF discretized \\ adjoint equations};
      \\
        \node [block_left] (assessment3) {Adams-Bashforth/Moulton\\ discretized\\ Control Problem}; 
      & \node [block_left] (excluded3) {Adams-Bashforth/Moulton\\ discretized\\  adjoint equations}; \\
       % enrollment - row 2
      };% end matrix

    \begin{scope}[every path/.style=line2]
            \path (referred) edge[above] node {} (excluded1);
    \end{scope}
    
     \begin{scope}[every path/.style=line3]
        \path (assessment) edge[above] node {(?)} (excluded2);
        \path (assessment3) edge[above] node {(?)} (excluded3);
    \end{scope}
    \begin{scope}[every path/.style=line]
            \path (referred) edge[above] node {}(assessment);
        \path (excluded1) edge[above] node {} (excluded2);
    \end{scope}

     \path (assessment)+(4.75,1.0) node {\bf Linear Multistep Methods} (excluded2)+(0.0,1.0);
    
% Draw background

  \begin{pgfonlayer}{background}
    % Left-top corner of the background rectangle
    \path (assessment)+(-2.75,0.85) node (assessment) {};
    % Right-bottom corner of the background rectangle
    \path (excluded3)+(+2.75,-0.9) node (excluded3) {};
    % Draw the background
    \path[fill=yellow!30,rounded corners, dashed,thick,draw=black!50,inner sep=6pt]
    (assessment) rectangle (excluded3) ;
  \end{pgfonlayer}

  \end{tikzpicture}
\caption{Time--dependent optimal control problems  discretized using linear--multi step methods.
Discretization of the arising adjoint equations either using discretized optimal control 
problems or discretized continuous adjoint equations \eqref{osOCPb}. We investigate the 
relation indicated by the question mark in the figure. }\label{overview}
\end{figure}
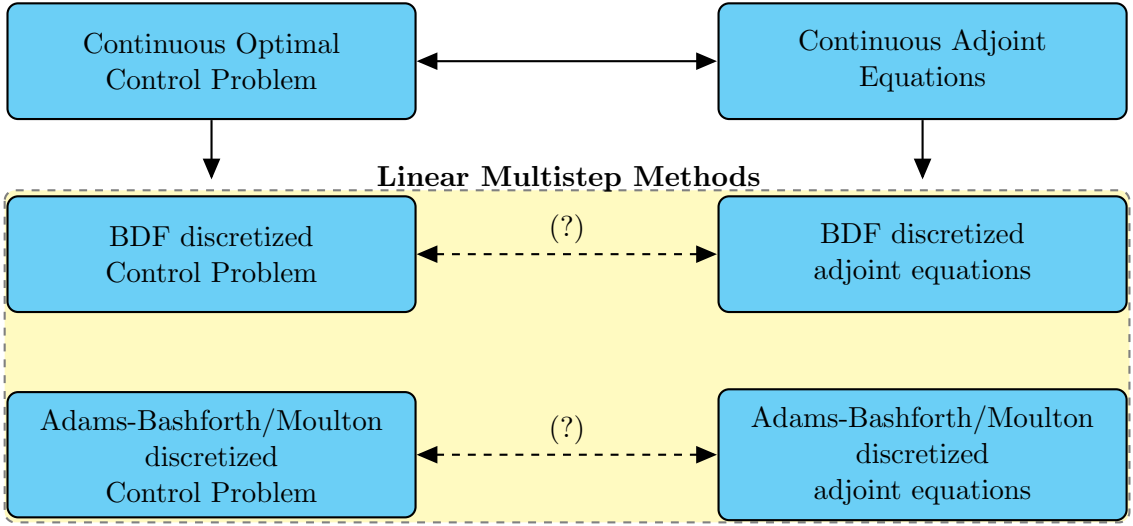

For possible numerical discretization we investigate the relations depicted in Figure \ref{overview}.
Therein,  we consider two different linear multi--step schemes for the discretization of 
the forward equation \eqref{osOCPa} and the adjoint equation \eqref{osOCPb}. Also, we consider
the optimality conditions \eqref{osOCPa}--\eqref{osOCPb} for the discretized problem. Then, we establish possible connections between both approaches. A similar investigation will be carried out 
for semi--Lagrangian discretization of hyperbolic relaxation systems. 

The ordinary differential equation is discretized using a linear multi--step method on $[0,T].$
For simplicity an equidistant grid in time $t_i = \dt \; i$ for $i=0,\dots,N$ such that $N \dt = T$
is chosen. The point value at the grid point $t_i$ is numerically approximated by $y(t_i)\approx y_i$,   $f(y(t_i)) \approx f(y_i)$ and $u(t_i) \approx u_i.$ A scheme is of order $p$ if the consistency error  of the numerical scheme is $y(t_i)=y_i + O(\dt^p),$ see \cite{HNW93}.  An $s-$stage linear multi-step scheme is defined by  \cite{HLW06,HNW93}  two vectors $a \in  \R^s, $ with components denoted by $a=(a_0,\dots, a_{s-1}),$ and $b \in \R^{s+1}$ with 
components $b=(b_{-1},b_0,\dots,b_{s-1}).$ Depending on the choice of $a,b$ we obtain so called Adams methods or BDF methods. In the case of BDF methods we have $b_i = 0, i\geq 0$ 
but $b_{-1}\not =0.$ Further, we define the 
numerical approximation of the solution at time $t_n,\dots,t_{n-s+1}$ as $Y_n=(y_n,\dots,y_{n-s+1}).$ For a $s-$stage multi-step scheme  we obtain an approximation to the solution 
$y(t)$ on the time interval $[ (1-s)\dt, \dots, T]$ that is denoted by
$ \vy=(y_{1-s},\dots,y_0,y_1,\dots,y_N). $
\par 

\subsection{Discretization of the optimal control problem}
The continuous problem \eqref{OCP} is discretized using an $s-$stage scheme. The initial condition is discretized by $Y_0=(y_0)_{i=0}^{s-1}$ and where 
$(y_0)_{s-1}$ for $i>0$ is an approximation to $y_0$. Further, $(y_0)_i$ for $i<s-1$ 
an approximation to the solution $y(t)$ of \eqref{osOCPa} at time $(1-s+i)\dt.$ In practice,
initialization may pose a difficulty and it can be observed that the order of scheme deteriorates
if the initialization has not been done properly. We assume a consistent  initialization at the order of the scheme.

Then, for a given control sequence $\vu:=(u_n)_n$ a linear multi--step discretization of equation \eqref{osOCPa} is of the following form 
\begin{equation}\label{scheme}
y_{n+1}= - a^t Y_n + \dt b^t F(Y_n,U_n), \; n \geq 0,
\end{equation}
where $F(Y_n,U_n) = ( f(y_i,u_i) )_{i=0}^{i=n-s+1}$. In order to compute the discretized linear multi--step optimality conditions it is advantageous to rewrite the 
previous system in matrix form 
\begin{equation}\label{matrix form}
\vy = - A \vy + \dt B F(\vy,\vu) + (Y_0,0,\dots,0)^t, 
\end{equation}
where 
 $A,B \in \mathbb{R}^{N+s \times N+s }$ have the same structure, namely, 
 \begin{align*} 
  A &= \begin{pmatrix}
 0_{s \times s} & \\
a_{s-1} & a_1 &  a_0 & 0_{s+1,s+1} \\
0 & a_{s-1} & a_1 &  a_0 & 0_{s+2,s+2} \\
 & \vdots \\
0 & \dots & 0 & a_{s-1} & a_1 &  a_0 & 0_{s+N,s+N} 
 \end{pmatrix}
\end{align*}
and
 \begin{align*}
 B &= \begin{pmatrix}
 0_{s \times s} & \\
b_{s-1} & b_1 &  b_0 & b_{-1} \\
0 & b_{s-1} & b_1 &  b_0 & b_{-1} \\
 & \vdots \\
0 & \dots & 0 & b_{s-1} & b_1 &  b_0 & b_{-1} 
 \end{pmatrix}.
\end{align*}
Finally, we discretize the cost functional $j.$  Several possibilities exist, the simplest one 
being $j(y(T)) \approx j(y_N).$ Other choices might include a polynomial reconstruction of
$j(y(T))$ using the $s-$stages $j(Y_N).$ We denote the numerical approximation 
of $j(y(T))$ by $j(Y_N).$ 

\begin{lem}\label{lem1}
Using an $s-$stage linear multi-step method  the discretized optimality system \eqref{OCP} with 
equi-distant temporal discretization $t_n=n \dt$ reads
\begin{equation}
\label{discOCP}
\min_{ \vy, \vu } j(Y_N) \mbox{ subject to } \eqref{scheme}, \; Y_0=(y_0)_{i=1-s}^{0}.
\end{equation}
The  discrete optimality conditions for $i=1-s,\dots,N$ are given by 
\begin{subequations}\label{discOCP}
\begin{align}
\vy &= - A \vy + \dt B F(\vy,\vu) + (Y_0,0,\dots,0)^t,  \\
0    &=  ( B^t \vp )_i f_u(y_i,u_i),   \\
0   &=  p_i + (A^t\vp)_i - \dt ( B^t \vp )_i f_y(y_i,u_i) + \partial_{y_i} j(Y_N). \label{discOCPc}
\end{align}
\end{subequations}
The initial conditions for $\vy$ are $y_i=(y_0)_i, i=1-s,\dots,0.$ The terminal condition
for multiplier $\vp$ are obtained from \eqref{discOCPc} for $i=N-s+1,\dots,N$ and read
e.g. for $i=N$
\begin{equation}\label{tdata}
 0 = p_N +  \partial_{y_N} j(Y_N) - b_{-1} \dt p_{N} f_y(y_N,u_N). 
 \end{equation}
\end{lem}
\noindent {\em Proof.} Due to the definition of a linear multi-step scheme the 
solution $\vy$ exists for any choice of $\vu.$ Therefore, we may write $\vy=\vy(\vu)$ 
and the constrained minimization problem \eqref{discOCP} reduces to an unconstrained
problem in $\vu.$ Hence, the discrete optimality conditions are necessary. They are derived
as saddle point of the discrete Lyapunov function given by 
$$L(\vy,\vu,\vp):=j(Y_N) + \vp^t \vy + (A^{t}\vp)^{t} \vy - \dt (B^{t} \vp)^{t} F(\vy,\vu) -  \vp^t (Y_0,0,\dots,0)^t, $$
where $\vp$ denotes the vector of adjoint states. Computing the partial derivatives of $L$
with respect to $\vu$ and $\vy$, respectively, yields the discrete optimality conditions where
we denote by $f_u$ and $f_y$ the partial derivatives of $f$ with respect to $u$ and $y.$ 
For the computation note that 
$$ (A^t \vp)_i  = a^t (p_{i+1},\dots,p_{i+s}), \mbox{ and }   (B^t \vp)_i = b^t  (p_i, p_{i+1},\dots,p_{i+s}).$$
Also note that the multipliers $p_i$ for $i=1-s,\dots,0$  only appear in the computation 
of $u_i$ for $i<0.$ Using the initial data $Y_0$ and the recalling the form of $A,$ we
observe that they do not enter the optimality conditions. Therefore, equation \eqref{discOCPc}
is in fact required only to hold for $i\geq0.$ $\hfill \blacksquare$
 
\begin{remark}\label{remark1}
It is important to remark that the equation \eqref{discOCPc} does in general 
{\em not} lead to  a linear multi--step method for the adjoint equation \eqref{osOCPb}. 
It utilized a fixed discretization point $f_y(y_i,u_i)$ even so $B^t p_i$ is the interpolation of
$p$ using values from $t_i,\dots, t_i+s \dt.$ 
\end{remark}

In view of Remark \ref{remark1} we consider a linear multi-step method applied to 
\eqref{osOCPb}. For notational simplicity we transpose equation \eqref{osOCPb} 
and obtain 
\begin{equation}
\label{cntadj}
 - p'(t) = f_y(y(t),u(t))p(t), \; p(T)=j_y(y(T)).
 \end{equation}

\begin{lem}\label{lem2}
A $s-$stage linear multi-step method applied to  equation \eqref{cntadj} on an equidistant 
grid for given functions $y(t),u(t)$ with discretizations $(\vy,\vu)$ is given by 
\begin{equation}\label{disccntadj} 
p_{n-1} = - \sum_{i=0}^{s-1} a_i  p_{n+i} + \dt b_i f_y(y_{n+i-1},u_{n+i-1})  p_{n+i}
\end{equation}
and terminal condition $P_N=\left( (j_y)(y_i)\right)_{i=N}^{N+s}.$
\end{lem}
\noindent {\em Proof.} We define $\g(t)=f_y(y(T-t),u(T-t))$ and $\p(t)=p(T-t)$ and obtain the equivalent equation
$$ \p'(t) = \g(t) \p(t), \; \p(0)=j_y(y(T)).$$

A linear multi-step method on the grid $t_i = i \; \dt$  for the adjoint variable
  $\p_n=\p(t_n)$ $\g_i=\g(t_i)$ is then given by
$$\p_{n+1} = -\sum_{i=0}^{s-1} a_i \p_{n-i} + \dt b_i \g_{n-i} \p_{n-i}$$
or transformed in original variables, i.e.,  $p_N=\p_1, p_1 = \p_N$, $ p_n = \p_{N-n+1}, \;  g_n = \g_{N-n+1}$,  reads as 

$$p_{n-1} = - \sum_{i=0}^{s-1} a_i  p_{n+i} + \dt b_i g_{n+i} p_{n+i} $$ 

Since $g_i = \g_{N-i+1}=f_y( y(T-t_{N-i+1}), u(T-t_{N-i+1}) ) = f_y( y_{i-1}, u_{i-1} )
$ we obtain the discretized continuous adjoint as \eqref{disccntadj}.
$\hfill \blacksquare$

Now, comparing \eqref{disccntadj} and \eqref{discOCPc} we observe that 
depending on $a$ $b,$ both equations are equivalent.

\begin{lem}\label{lem3}
Assume $j(y(T))$ is approximated by $j(y_N).$ Then, for $t<T$, the update formula for    discretize--then--optimize, i.e.,
equation \eqref{discOCPc} and optimize--then--discretize \eqref{disccntadj}
 coincide up to $O(\Delta t^p)$ for  BDF type methods.
\end{lem}
\noindent {\em Proof} In case of BDF methods we have $b_i=0$ for $i\geq 0.$ Therefore, 
 equation \eqref{osOCPc} reads for $i<N$:
$$p_{n-1} = - \sum_{i=0}^{s-1} a_i  p_{n+i} + \dt b_i f_y(y_{n+i-1},u_{n+i-1})  p_{n+i} $$ 
On the other hand, \eqref{disccntadj} reads 
$$p_{n-1} = - \sum_{i=0}^{s-1} a_i  p_{n+i} + \dt b_{-1} f_y(y_{n-2},u_{n-2})  p_{n-1} $$ 
Since $y_{n-2} = y_{n-1} + O(\dt^p)$ the equations
 coincide up to the order of the scheme for $i<N.$   $\hfill \blacksquare$
 
\begin{remark}
%We list some considerations about the proposed schemes
% \begin{enumerate}
% \item 
The terminal  data is discretized in the case of Lemma \ref{lem1} by \eqref{tdata} and by $p_N=\partial_{y_N} j(y_N)$ in  the case of Lemma \ref{lem2}. However, for the continuous discretization of the adjoint equation \eqref{lem2}  this choice can be altered  to be consistent with the discretization of Lemma~\ref{lem1}. 
%  \item 
Clearly, if $f_y=const,$ different discretizations do not affect the method. Therefore, the previous Lemma only states necessary conditions. We refer to Section \ref{resultsycst} for numerical results. 
%  \item 

We further observe that {\em no} method with $b_i \not = 0$ for $i\geq 0$ yields a consistent discretization in both approaches. Hence, in Figure \ref{overview} only the question mark in between the BDF methods can be answered positive.  In fact, for Adams--Bashfort and Adams--Moulton type   methods we observe a decay in the order, see Section \ref{decayMoulton}.
%\item
The results presented in \cite{Sandu} also show that in general one can only expect
first--order convergence without further assumptions on the choices of $a$ and $b.$

%\item
Finally, in \cite{HertyPareschiSteffensen2013} also the question of long--term integration
of the optimality conditions has been studied. In the context of linear multi--step scheme
it is already known \cite{HNW93} that there is {\em no } high--order scheme that is symplectic. 

%\end{enumerate}
 \end{remark}

%%%%%%%%%%%%%%%%%%%%%%%%%%%%%%%%%%%%%%%%%%%%%%%%%%

\section{Linear multi-step methods for optimal control problems of relaxation systems}

\subsection{Semi-lagrangian schemes for relaxation approximations}
Relaxation approximations to hyperbolic conservation laws have been introduced in \cite{JX}. To exemplify the approach we consider a nonlinear
scalar conservation law of the type 
\begin{equation}\label{cons law}
u_t + F(u)_x = 0, \quad x \in \R, t \geq 0
\end{equation} 
and initial datum $u(0,x)=u_0(x).$ The flux function $F:\R \to \R$ is assumed to 
be smooth. In order to apply a numerical integration scheme we introduce a relaxation approximation 
to \eqref{cons law} as 
\begin{equation}
\begin{split}
\label{JX} u_t + v_x &= 0, \\ v_t + a^2 u_x &= \frac{1}\epsilon \left( F(u)-v \right).
\end{split}
\end{equation}
Note that the above approximation can be interpreted as a BGK-type kinetic model \cite{2} by introducing the Maxwellian 
equilibrium states $E_{f}$ and $E_g$ given by 
$$ E_{f}(u) = \frac{1}{2a} \left( au + F(u) \right), \; E_g(u) = \frac{1}{2a} \left(au- F(u) \right).$$
The kinetic variables $f,g:\R^+ \times \R \to \R$ fulfill then 
\begin{equation}
\begin{split}
\label{bgk} f_t + a f_x &= \frac{1}\epsilon \left( E_f(u) -f \right), \\ g_t - a g_x &= \frac{1}\epsilon \left( E_g(u) - g \right),
\end{split}
\end{equation}
with $u = f + g $ and $v=a(f-g)$. Herein, $a$ is the characteristic speed of the transported variables and it is assumed that this speed bounds the 
eigenvalues of \eqref{cons law}, i.e., the subcharacteristic condition holds 
$$ a \geq \max_{ x \in  \R } | F'(u_0(x)) |.$$

In the formal relaxation limit $\epsilon \to 0$ we recover the following relations
\begin{equation}
\label{bgk limit} 
f=E_f(u),\quad g=E_g(u), \quad v=a(f-g)=F(u).
\end{equation}
Therefore, $u=f+ g$ fulfills in the small relaxation limit the conservation law \eqref{cons law}.
Due to the linear transport structure in equation \eqref{bgk} semi--Lagrangian schemes can be used and the system
\eqref{bgk} reduces formally to a coupled system of ordinary differential equations. Let us mention that recently, linear multi-step methods have been proposed to numerically solve kinetic equations of BGK-type \cite{GroppiRussoStracquadanio}. 

Let 
$$ \bar{f}(t,y):= f(t,y+at),\quad \bar{g}(t,y) = g(t,y-at)$$
for a point $y\in \R.$ Then, the macroscopic variable $u(t,x)$ is obtained through 
$$
u(t,x) = \bar{f}(t, x-at) + \bar{g}(t,x+at),
$$
and for any $y$ we have $$ \frac{d}{dt} \bar{f}(t,y) = f_t(t,y+at) + a f_y(t,y+at).$$
Therefore, the unknowns $\bar{f}$ and $\bar{g}$ fulfill a coupled system of ordinary differential equations for all $y \in \R:$
\begin{subequations}\label{bgk ode}
\begin{eqnarray}
\frac{d}{dt} \bar{f}(t,y) = \frac{1}\epsilon \left( E_f(u(t,y+at)) - \bar{f}(t,y) \right),\quad  u(t,y+at) = \bar{f}(t,y) + \bar{g}(t,y+2at)  \\ 
\frac{d}{dt} \bar{g}(t,y) = \frac{1}\epsilon \left( E_g(u(t,y-at)) - \bar{g}(t,y) \right),\quad u(t,y-at) =  \bar{f}(t,y-2at) + \bar{g}(t,y).
\end{eqnarray}
\end{subequations}

Next, we turn to the numerical discretization of the previous system of (parameterized) ordinary differential equations. We introduce
a spatial grid of width $\dy$ and denote for $i \in \mathbb{Z}$ the grid point $y_i = i \dy.$ Similarly, in time we introduce
a spatial grid of width $\dt$ and denote by $t_n= n \dt$ for $ n \in \mathbb{N}.$ 

Note that explicit schemes require a CFL condition for the relation between 
spatial and temporal grid to hold, i.e., 
\begin{equation} \dt \leq a \dy \label{CFL}. \end{equation}
In the case of implicit discretizations as e.g. BDF this  is not required.  The point values of $\bar{f}$ and $\bar{g}$ are denoted by 
\begin{equation*}
 \f^n_i = \f(t_n,y_i) := \f( n \dt, i \dy), \quad \g^n_i =\g(t_n,y_i) := \g(n\dt,i \dy).
\end{equation*}
For each $y_i$ we apply a linear--multi step scheme to discretize in time.
For simplicity here we restrict the analysis to BDF methods. These require only a single
evaluation of the source term and this evaluation is implicit. Therefore, the time discretization $\dt$ does {\em not }
dependent on the size of $\epsilon.$  For an $s-$stage scheme and using a temporal discretization  $\dt=a\dy$ we obtain an explicit scheme on the indices $i,$ given by  
\begin{subequations}\label{disc01}
\begin{eqnarray}
\f^{n+1}_i = \frac{\dt b_{-1} }{\dt b_{-1} +\epsilon}  E_f\left( \f^{n+1}_i + \g^{n+1}_{i+2(n+1)}\right)  - \frac{\epsilon}{\dt b_{-1} +\epsilon} \sum\limits_ {\ell= 0}^{s-1} a_{\ell} \f^{n-\ell}_i 
, \\
\g^{n+1}_i = \frac{\dt b_{-1} }{\dt b_{-1} +\epsilon}  E_g\left( \f^{n+1}_{i-2(n+1)} + \g^{n+1}_{i} \right) - \frac{\epsilon}{\dt b_{-1} + \dt} \sum\limits_ {\ell= 0}^{s-1} a_{\ell} \g^{n-\ell}_i.
\end{eqnarray}
\end{subequations}
Since there is no spatial reconstruction it suffers in the case of strong discontinuities in the spatial variable as observed in 
\cite{GroppiRussoStracquadanio}. 
\par 
We further investigate the  continuous system \eqref{bgk ode} 
and its discretization \eqref{disc01} in the particular case 
$$F(u) = c \; u, \; c>0.$$
For the relaxation system to approximate the conservation law we 
require $a \geq c.$ Using the semi-Lagrange scheme we observe that
 the choice $a=c$  leads to an exact scheme. In this  case we obtain 
  $E_f(u)=u$ and $E_g(u)=0.$ Furthermore, the equations \eqref{bgk ode}
  reduce to 
  \begin{equation}
  \begin{split}
  \frac{d}{dt} \bar{f}(t,y) &= \frac{1}\epsilon \bar{g}(t,y+2at),\\
  \frac{d}{dt} \bar{g}(t,y)&=-\frac{1}\epsilon \bar{g}(t,y).
  \end{split}
  \end{equation}
As initial data for $\bar{f}$ and $\bar{g}$ we may chose $\bar{f}(0,x)=u_0(x)$ and $ \bar{g}(0,x)=0$. Then, the previous dynamics yield in the limit $\epsilon \to 0$ the projections $\bar{g}(t,y)=0$ and 
$\bar{f}(t,y)=u_0(y)$. Rewritten in Eulerian coordinates we obtain 
$ u(t,x) = u(t,x-a t)$ being the solution to the original linear transport equation \eqref{cons law} if $a=c$. This computation shows that $a=c$ is necessary for consistency with the original problem in the small $\epsilon$ limit.  The discretized equations \eqref{disc01}  with initial data $\g^0_i=0, \f^0_i = u_0(x_i)$ simplify to 
$\g^n_i\equiv 0$ and 
\begin{equation}\label{disc f}
\f^{n+1}_i \left( 1 - \frac{\dt}{ \dt b_{-1} + \epsilon} \right) = - \frac{\epsilon}{ \dt b_{-1} + \epsilon} \sum\limits_{\ell=0}^{s-1} a_\ell \f^{n-\ell}_i.
\end{equation}

Summarizing, equation \eqref{disc f} shows that the BDF discretization in the case of a linear 
flux function with suitable initialization of the relaxation variables leads to 
a high--order formulation in Lagrangian coordinates. The discretization is
independent of the spatial discretization and there is {\bf no} CFL condition. 

However, this discretization is only exact in the case of a linear transport equation. In the case $F(u)$ nonlinear additional interpolation needs to 
be employed. Then, due to the Lagrangian nature of the scheme, the spatial
resolution and the temporal is coupled through the interpolation. 
%Therefore, we only present qualitative results in the nonlinear case.

%{\bf Lorenzo: I have some doubts on the next section. I suggest we consider the case of hyperbolic relaxation systems of dimension $N$ with a linear (BGK-type) relaxation term and without discussing the limit $\varepsilon \to 0$ which is not the goal of our manuscript. Note that we should have one test case with $N > 2$ for example Broadwell with $N=3$ is very close to Burgers (I can discuss this with Giacomo)}
%
%{\bf Michael: I changed everything to two velocties. the limit eps to zero is ok, since the adjoint equation is linear.   }

\subsection{Derivation of adjoint equations for the control problem} 

We will derive the {\em adjoint } BDF schemes for the previous discretization and we compare the discrete adjoint equations with the formal continuous adjoint equation to the conservation law \eqref{cons law}. In order to simplify notations, we denote the spatial
variable in the kinetic and Lagrangian frame also by $x$ (instead of $y$). Furthermore, in view of generalizations to the case of systems with a larger number of velocities, we introduce the velocities $v_1=a, v_2=-a$ 
as well as the kinetic variables $f^1=f$ and $f^2=g$ and the corresponding equilibrium as $E_1=E_f$ and $E_2=E_g.$

% We consider a high--order BDF forward scheme for  BGK approximation 
% to \ref{cons law} with velocities $v_1=a$ and $v_2=-a.$ The corresponding kinetic
%variables are $f^j=f^j(t,x)$, initial data $f^j_0$, and the 
%Maxwellian $E_j$ and the macroscopic quantity $u(t,x)$. 

Then, the  hyperbolic relaxation approximation is  given by the kinetic transport equation for $j=1,2$
\begin{subequations}\label{BGKforward}
\begin{eqnarray}
f_t^j + v_j f^j_x &=& \frac{1}\epsilon \left( E_j(u) - f^j \right), 
\\
f^{j}(0,x) &=& f^{j}_0(x),
\end{eqnarray}
\end{subequations}
with $u(t,x) = \sum_{j} f^{j}(t,x)$. We recall that the local equilibrium states have the property $\sum_j E_{j}(u)=u$ that will be used in the differential calculus later on.

As before, we define the Lagrangian variables $\fl$ as 
$\fl^j(t,x) = f(t,x+v_j t)$ and the macroscopic quantity $u$ as 
$u(t,x) = \sum_{j} \fl^j(t,x-v_jt).$
Then, equation \eqref{BGKforward} is equivalent to the 
ODE system \eqref{ODE} and initial data $\fl^j(0,x)=f^j_0(x).$ 
\begin{eqnarray}\label{ODE}
\partial_t \fl^j(t,x) = \frac{1}\epsilon \left( E_j(u(t,x+v_j t)) - \fl^j(t,x) \right).
\end{eqnarray}
Consider the integral form of \eqref{ODE} on the time interval $[s,t].$ Since $f^j(t,x)=\fl^j(t,x-v_j t)$ 
we have  for all $s<t$ and all $x \in \R$:
\begin{eqnarray*}
% \fl^j(t,x - v_j t ) -\fl^j(s, x - v_j t ) = \frac{1}\epsilon \int_s^t  E_j(u(\tau,x - v_j t+v_j \tau)) - \fl^j(\tau,x -  v_j t) ds \\
 f^j(t,x ) - f^j(s, x- v_j (t-s) ) = \frac{1}\epsilon \int_s^t  
 E_j(u(\tau,x - v_j ( t- \tau )) -  f^j(\tau,x - v_j (t-\tau) ) d\tau
\end{eqnarray*}

Upon summation on $j$ we have 
%\begin{eqnarray*}
%\frac{1}N\sum_j \left(  f^j(t,x ) - f^j(s, x- v_j (t-s) ) \right)= \\
%\frac{1}\epsilon \int_{s}^t   \frac{1}N \sum_j \left(
% E_j(u(\tau,x - v_j ( t- \tau )) -  f^j(\tau,x - v_j (t-\tau) ) \right) d\tau =0.
%\end{eqnarray*}
%Therefore, 
for $s<t$ we have
$$u(t,x) = \sum_j f^j(s, x- v_j (t-s)).$$
%Hence, we may evaluate the Maxwellian at time $t$ using data 
%at previous time. It is therefore not necessary to solve a nonlinear system within a BDF--type method.

We are interested in initial conditions $f^j_0(\cdot)$ minimizing a cost function $J$ 
depending on the macroscopic variables $u_0 = \frac{1}N\sum_j f^j_0(x)$ 
as well as $u(T,x)$ at some given point $T>0.$ The dynamics of $u$ is approximated by the BGK formulation \eqref{BGKforward}. 
\begin{equation}\label{OPTx}
\min\limits_{ f^{j}_0(x), j=1,2} 
\int J(u(T,x),u_0) dx\mbox{ subject to } \eqref{BGKforward}.
\end{equation}
It is straightforward to derive the formal optimality conditions including the 
 formal adjoint equations for the variables $\lambda^j(t,x).$ Those are defined up to 
 a constant and therefore we state the adjoint equation in the re-scaled variables
 $\frac{1}2 \lambda^j(t,x)$ for $j=1,\ldots,N$ as follows 
 \begin{equation}
\begin{aligned}\label{BGK adjoint}
- & \lambda^{j}_t - v_j   \lambda^{j}_x = - \frac{1}\epsilon \left( \lambda^{j} - \sum_k \lambda^k  E_k'(u(t,x)) \right), \\
 &\lambda^{j}(T,x)  + J_u(u(T,x),u_0) = 0.
 \end{aligned}
 \end{equation}
 The adjoint multipliers and the optimal control $u_0$ are then related according
 to 
$$ -   \lambda^j(0,x) + J_{u_0}(u(T,x),u_0) = 0.$$
The property of the local equilibrium implies 
$ \sum_j E_j'(u) = 1$ and therefore, 
$$ \sum_j E_j'(u) \left(  - \lambda^{j}_t - v_j \lambda^{j}_x \right) = 
0.$$
In the formal limit  $\epsilon \to 0$ we obtain that $\lambda^{j} = \sum_k \lambda^{k} E_k'(u) =  \lambda$
and therefore $\lambda^j$ is independent of $j.$  

\begin{lem}
Up to $O(\epsilon^2)$ the equations \eqref{BGK adjoint} are a viscous 
approximation to the linearized adjoint equation to equation \eqref{cons law} given by
$$ - p_t - F'(u)p_x=0.$$
\end{lem}
\noindent {\bf Proof.} 
For the local equilibrium $E_{j}$ it holds $E_1(u) + E_2(u) = u$ and additionally $ E_1(u) - E_2(u) = F(u)/a$, for all $u\in \R$, and therefore, $ E_1'(u) - E_2'(u) = F'(u)/a.$ We denote by $\lambda^\pm=\lambda^{1,2}.$ We obtain for the sum and the difference of $\lambda^\pm \pm \lambda^\mp$ 
the following equations 
\begin{eqnarray*}
- (\lambda^++\lambda^-)_t - a (\lambda^+-\lambda^-)_x =  \frac{1}\epsilon \left( F'(u) /a \right) (  \lambda^+ - \lambda^- ), \\
- (\lambda^+ - \lambda^-)_t - a (\lambda^+ + \lambda^-)_x = -\frac{1}\epsilon (\lambda^+ - \lambda^-). 
\end{eqnarray*}
Denote by $ \lambda = \lambda^++\lambda^-$ and by $\phi:=\lambda^+-\lambda^-.$ Then, the equations are
equivalent to 
$$ -\lambda_t - a \phi_x =  \frac{1}{\epsilon a} F'(u) \phi,  \; - \phi_t - a \lambda_x = - \frac{1}\epsilon \phi.$$
Hence,  $ \phi =  \epsilon ( a \lambda_x) + O(\epsilon^2)$ and therefore,  $ - \lambda_t - F'(u) \lambda_x  = \epsilon a^2 \lambda_{xx}. $ $\hfill \blacksquare$

Next, we discuss BDF discretization of the adjoint equations. The adjoint variables
$\lambda^j$ are transported backwards in space and time. In order to derive a semi--Lagrangian description we define
$$\la^{j}(t,x) = \lambda(t, x+v_j t)$$
and define the terminal  data as $\la^{j}(T,x) = - J_u(  u(T,x+v_j T), u_0(x+v_j T) ).$
The semi--Lagrangian formulation of the adjoint equation is
\begin{equation} \label{adjoint}
 - \partial_t \la^{j}(t,x) =  - 
\frac{1}\epsilon
\left(  \la^j(t,x) - \sum_k \lambda^k(t,x+v_j t) 
E_k'(u(t,x+v_j t)) \right),
\end{equation}
or upon integration from $s$ to $t$ with $s<t$
$$ \la^{j}(s,x) - \la^{j}(t,x) =  -\frac{1}\epsilon
\int_s^t \left(  \la^j(\tau,x) - \sum_k \lambda^k(\tau,x+v_j \tau) 
E_j'(u(\tau,x+v_j \tau)) \right) d\tau.$$
A BDF integrator with $s-$stages applied to this equation yields the discretized 
equation 
\begin{align}\label{eq adj lambda}
\la^j(t_{n-1},x) + \sum_{i=0}^{s-1} a_i \la^j(t_{n+i},x) =- \frac{\dt b_{-1}}\epsilon \left( \la^j(t_{n-1},x) - 
Z(t_{n-1},x+v_j t_{n-1}) \right)
\end{align}
where the source term is given by 
$$Z(t,y) := \sum_k \lambda^{k}(t,y) E_k'(u(t,y)).$$

Similarly to the forward equations we evaluate $Z$ without  knowledge on $\lambda^{k}(t,y)$ using the integral formulation of the problem above.
We show this relation in the time--discrete case. Denote the 
discretize Eulerian adjoint variables by 
$ \la^j(t_{n+i},x-v_j t_{n-1}) = \lambda^{j}(t_{n+i}, x+v_j t_{n+i} - v_j t_{n-1})
$
where $t_{n+i}=t_{n-1}+(i+1)\dt, \; i=0,1\dots,s-1.$ Then, 
\begin{align*}
\lambda^j(t_{n-1},x) + \sum_{i=0}^{s-1} a_i \lambda^j(t_{n+i},x+v_j(i+1) \dt) =- \frac{\dt b_{-1}}\epsilon \left( \lambda^j(t_{n-1},x) -   Z(t_{n-1},x) \right).
\end{align*}
After multiplication with $E_j'(u)$ and summation on $j$ we obtain
\begin{align*}
Z(t_{n-1},x) + \sum_j \sum_{i=0}^{s-1}  E_j'(u(t_{n-1},x))a_i \lambda^j(t_{n+i},x+v_j(i+1) \dt) = \\
- \frac{\dt b_{-1}}\epsilon \left( Z(t_{n-1},x) -   \sum_j E_j'(u(t_{n-1},x))  Z(t_{n-1},x) \right) = 0.
\end{align*}
%Z(t_{n-1},x) = 
%- \frac{1}N \sum_j \sum_{i=0}^{s-1}  E_j'(u(t_{n-1},x) )a_i \lambda^j(t_{n+i},x+v_j(i+1) \dt) 
%\end{align*}
The  equation for $\lambda^{j}(t_{n-1},x)$ 
is explicit since  $Z(t_{n-1},x)$  depends only on $\lambda^j(t_{n+i},\cdot)$ for $i\geq 0.$  Equation \eqref{eq adj lambda} is equivalent to 
\begin{align*}
\lambda^j(t_{n-1},x)   \frac{\epsilon + \dt b_{-1}}\epsilon  
=-  \sum_{i=0}^{s-1} a_i \lambda^j(t_{n+i},x+v_j(i+1) \dt) + 
\frac{\dt b_{-1}}{\epsilon} 
Z(t_{n-1},x), 
\end{align*}
where 
$$\frac{\dt b_{-1}}{\epsilon} Z(t_{n-1},x) = - \frac{\dt b_{-1}}{\epsilon} 
  \sum_j \sum_{i=0}^{s-1}  E_j'(u(t_{n-1},x) )a_i \lambda^j(t_{n+i},x+v_j(i+1) \dt).$$
  Therefore the  adjoint BDF discretization of the continuous adjoint equations
  in Eulerian coordinates is given by 
\begin{subequations}\label{implement adjoint}
 \begin{align}
  \lambda^j(t_{n-1},x) = - \frac{ \epsilon }{ \epsilon + \dt b_{-1} }
  \sum_{i=0}^{s-1} a_i \lambda^j(t_{n+i},x+v_j(i+1) \dt) -  \\
  \frac{ \dt b_{-1}}{ \epsilon + \dt b_{-1} } \sum_j \sum_{i=0}^{s-1}  E_j'(u(t_{n-1},x) )a_i \lambda^j(t_{n+i},x+v_j(i+1) \dt). 
\end{align}
\end{subequations}
We observe that the limit $\epsilon\to 0$ exists and it is independent of $\lambda^j$ 
as in the continuous case. Further, for $\epsilon>0$ and $\dt \to 0$ we obtain 
the interpolation property of BDF methods, i.e.,   
\begin{align*} 
\lambda^j(t_{n-1},x)    
=-   \sum_{i=0}^{s-1} a_i \lambda^j(t_{n+i},x). 
 \end{align*}
 
Summarizing, the adjoint equation \eqref{BGK adjoint} can be solved
efficiently using any BDF scheme in the formulation \eqref{implement adjoint}. 

\begin{lem}  
Consider the the adjoint equation \eqref{BGK adjoint} for the unknown 
adjoint variables $\lambda^1$ and $\lambda^2.$   Then, the scheme given by  \eqref{implement adjoint} is a discretization of the 
adjoint equation using a linear multi--step scheme of the family of BDF schemes.  In the limit $\Delta t \to 0$ and for $\epsilon>0$ this discreitzation is
 consistent with the interpolation property of BDF schemes.
\end{lem}

\subsection{Generalization to systems of conservation laws}
The approach here described can be extended to general one-dimensional hyperbolic relaxation systems and kinetic equations of the form \cite{2, JX}
\begin{subequations}\label{BGKN}
\begin{eqnarray}
f_t^j + v_j f^j_x &=& \frac{1}\epsilon \left( E_j({\bf u}) - f^j \right),\quad j=1,\ldots, N 
\\
f^{j}(0,x) &=& f^{j}_0(x),
\end{eqnarray}
\end{subequations}
where now ${\bf u}$ is a $n$-dimensional vector with $n < N$, such that there exists a constant matrix $Q$ of dimension $n\times N$ and ${\rm Rank}(Q)=n$ which gives $n$ independent conserved quantities ${\bf u}=Q{\bf f}$, ${\bf f}=(f^1,\ldots,f^N)^T$. Moreover, we assume that there exist a unique local equilibrium vector such that $Q{\bf E(u)}={\bf u}$, ${\bf E(u)}=(E_1({\bf u}),\ldots,E_N({\bf u}))^T$. 

From the properties of $Q$, using vector notations, we obtain a system of conservation laws which is satisfied by every solution of \eqref{BGKN} 
\begin{equation}
Q{\bf f}_t + QV{\bf f}_x = 0,
\end{equation}
where ${V}={\rm diag}\{v_1,\ldots,v_N\}$.
For vanishing values of the relaxation parameter $\varepsilon$ we have ${\bf f=E(u)}$ and system \eqref{BGKN} is well approximated by the closed equilibrium system
\begin{equation}
{\bf u}_t + F({\bf u})_x = 0,
\end{equation}
with $F({\bf u})=QV{\bf E(u)}$.
Using these notations, the control problem detailed in this Section corresponds to $N=2$, $n=1$ and $Q=(1,1)$.

%%%%%%%%%%%%%%%%%%%%%%%%%%%%%%%%%%%%%%%%%%%%%%%%%%
\section{Numerical results} 
We prove numerically previous results for BDF, Adams--Bashforth/Moulton integrators, for ODEs systems and relaxation  systems, presenting order of convergence and qualitatively results. 
We refer to Appendix \ref{tables} for a detailed definition of BDF, Adams--Bashforth/Moulton integrators.

\subsection{Convergence order for BDF and Adams--Bashfort/Moulton integrators}  \label{resultsycst}
% code: ODE/testsimple
In this section we verify the implementation of BDF and Adams--Bashfort/Moulton integrators
for the adjoint equation \eqref{osOCPb}. As discussed in Lemma \ref{lem1} to Lemma \ref{lem3} 
the derived adjoint schemes might be different depending on the approach taken in Figure \ref{overview}. However, in the special case $f_y=cst$ both approaches yield the same discretization
scheme and we do not expect any loss in the order of approximation. To illustrate we consider
$f_y=1$ and terminal data $p(T)=0.$ Then, the exact solution to equation \eqref{osOCPb} is given by 
$$ p(t)=\exp((T-t)).$$ The error is measured with respect to the exact solution. The results
are given in Table \ref{tab1}. The expected convergence order is
numerically observed for all tested methods. We only show the Adams--Bashfort and Adams--Moulton simulations. 

 \begin{table}[t]
 \centering
 \caption{Number of discretization points in time $N$, error in $L^\infty(0,T)$ for the approach discretize--then--optimize (Lemma \ref{lem1}) is shown in $L^\infty p$ with 
corresponding rate (Rate) and error in $L^\infty(0,T)$ for the approach optimize--then--discretize (Lemma \ref{lem2}) is shown in $L^\infty p$ with 
corresponding rate (Rate). We report from top to bottom different schemes: Explicit Euler, Adams--Bashforth(3), and Adams-Moulton(4).} \label{tab1}
\vspace{+0.25cm}
 \begin{tabular}{c|l| rr|rr}
 \hline
~ & \multicolumn{1}{c|}{$N$} & \multicolumn{1}{c}{$L^\infty p$} & \multicolumn{1}{c|}{Rate} & \multicolumn{1}{c}{$L^\infty  p(t)$} & \multicolumn{1}{c}{ Rate}\\
 \hline
 \hline
 \parbox[t]{4.5mm}{\multirow{5}{*}{\rotatebox[]{90}{\parbox[c]{1.6cm}{\footnotesize Explicit--Euler}}}} &40 &  0.0203478 & 2.11057 &  0.0203478 & 2.11057 \\
& 80 &  0.00490164 & 2.05354 &  0.00490164 & 2.05354 \\
 &160 & 0.00120324 & 2.02634 &  0.00120324 & 2.02634 \\
 &320 & 0.000298097 & 2.01307 &  0.000298097 & 2.01307 \\
 &640 &  7.41889e-05 & 2.00651 &  7.41889e-05 & 2.00651 \\
  \hline
% \end{tabular}
%\\
%~
%\\
%  \begin{tabular}{c|l| rr|rr}
 \hline
~ & \multicolumn{1}{c|}{$N$} & \multicolumn{1}{c}{$L^\infty p$} & \multicolumn{1}{c|}{Rate} & \multicolumn{1}{c}{$L^\infty  p(t)$} & \multicolumn{1}{c}{ Rate}\\
 \hline
 \hline
  \parbox[t]{4.5mm}{\multirow{5}{*}{\rotatebox[]{90}{\parbox[c]{2.cm}{\footnotesize Adams--Bashforth(3)}}}}&40 &  9.46513e-05 & 4.24563 & 9.46513e-05 & 4.24563 \\
 &80 &  5.42931e-06 & 4.12378 &  5.42931e-06 & 4.12378 \\
 &160 & 3.25127e-07 & 4.06169 &  3.25127e-07 & 4.06169 \\
 &320 & 1.9892e-08 & 4.03074 &  1.9892e-08 & 4.03074 \\
 &640 & 1.2301e-09 & 4.01534 &  1.2301e-09 & 4.01534 \\
 \hline
% \end{tabular}
%\\
%~
%\\
%   \begin{tabular}{c|l| rr|rr}
 \hline
~ & \multicolumn{1}{c|}{$N$} & \multicolumn{1}{c}{$L^\infty p$} & \multicolumn{1}{c|}{Rate} & \multicolumn{1}{c}{$L^\infty  p(t)$} & \multicolumn{1}{c}{ Rate}\\
 \hline
 \hline
 \parbox[t]{4.5mm}{\multirow{5}{*}{\rotatebox[]{90}{\parbox[c]{2.cm}{\footnotesize Adams--Moulton(4)}}}}&40 &  2.91401e-08 & 6.39089 &  2.91401e-08 & 6.39089 \\
 ~&80 &  3.99048e-10 & 6.1903 &  3.99048e-10 & 6.1903 \\
 ~&160 &  5.84258e-12 & 6.09381 &  5.84258e-12 & 6.09381 \\
 ~&320 &  8.86503e-14 & 6.04234 &  8.86503e-14 & 6.04234 \\
 ~&640 &  1.41997e-15 & 5.96419 &  1.41997e-15 & 5.96419 \\
 \hline
 \end{tabular}
 \end{table}

%\clearpage

\subsection{Loss of convergence order for Adams--Moulton integrators}  \label{decayMoulton}
% code: ODE/main

Compared to \eqref{resultsycst} we modify the adjoint equation by assuming 
$$f_y(y,u) = y(t),\qquad y(t)=t^2.$$
Terminal data for $p$ is again $p(T)=0.$ The exact solution of the adjoint equation is
explicitly known in this case and given by $p(t)=\exp((T^3-t^3)/3).$ Errors are measured with respect to the exact solution.  In view of Lemma \ref{lem3} we
expect only the BDF scheme to retain the high--order. The Adams--Moulton integrators
have $b_i\not =0$ for $i\geq 0$ and therefore the approach discretize--then--optimize
leads to inconsistent discretization of the adjoint equation \eqref{osOCPb}, see Lemma \ref{lem1}.
We show three different schemes: an explicit Euler, BDF(4) and Adams--Moulton(4). 
For each scheme we implement both versions, i.e., discretize--then--optimize and 
optimize--then--discretize. Clearly, in the case of the BDF method there is no difference
as expected due to Lemma \ref{lem3}. Also, for first--order methods there is no difference  
since $b_0=0.$ However, for the Adams--Moulton method  we observe the decay in 
approximation order in the case discretize--then--optimize. 
The results are given in Table \ref{tab2}. 
Obviously, we expect the same decay for Adams--Bashfort formulas. Those numerical results are skipped for brevity.

\begin{table}[tb]
 \centering
 \caption{Number of discretization points in time $N$, error in $L^\infty(0,T)$ for the approach discretize--then--optimize (Lemma \ref{lem1}) is shown in $L^\infty p$ with 
corresponding rate (Rate) and error in $L^\infty(0,T)$ for the approach optimize--then--discretize (Lemma \ref{lem2}) is shown in $L^\infty p$ with 
corresponding rate (Rate). We report from top to bottom different schemes: Explicit Euler, BDF(4), and Adams-Moulton(4).} \label{tab2}
\vspace{+0.25cm}
 \begin{tabular}{c|l| rr|rr}
 \hline
~ & \multicolumn{1}{c|}{$N$} & \multicolumn{1}{c}{$L^\infty p$} & \multicolumn{1}{c|}{Rate} & \multicolumn{1}{c}{$L^\infty  p(t)$} & \multicolumn{1}{c}{ Rate}\\
 \hline
 \hline
 \parbox[t]{4.5mm}{\multirow{5}{*}{\rotatebox[]{90}{\parbox[c]{1.6cm}{\footnotesize Explicit--Euler}}}}& 40 & 0.0358346 & 2.12096 &  0.00497446 & 1.76334 \\
 &80 &  0.00856002 & 2.06567 &  0.00144543 & 1.78305 \\
& 160 &   0.00209021 & 2.03397 &  0.000380452 & 1.92571 \\
 &320 &  0.00051634 & 2.01725 &  9.71555e-05 & 1.96935 \\
 &640 &  0.00012831 & 2.00869 &  2.45239e-05 & 1.98611 \\  \hline
% \end{tabular}
%\\
%~
%\\
%  \begin{tabular}{c|l| rr|rr}
 \hline
~ & \multicolumn{1}{c|}{$N$} & \multicolumn{1}{c}{$L^\infty p$} & \multicolumn{1}{c|}{Rate} & \multicolumn{1}{c}{$L^\infty  p(t)$} & \multicolumn{1}{c}{ Rate}\\
 \hline
 \hline
\parbox[t]{4.5mm}{\multirow{5}{*}{\rotatebox[]{90}{\parbox[c]{2.cm}{\footnotesize BDF(4)}}}}&
40 &  4.79238e-05 & 4.74597 &  4.79238e-05 & 4.74597 \\
 &80  &  1.35856e-06 & 5.1406 &  1.35856e-06 & 5.1406 \\
 &160  &  3.90305e-08 & 5.12133 &  3.90305e-08 & 5.12133 \\
 &320  &  1.16026e-09 & 5.07209 &  1.16026e-09 & 5.07209 \\
 &640  &  3.52961e-11 & 5.03879 &  3.52961e-11 & 5.03879 \\
   \hline
% \end{tabular}
%\\
%~
%\\
%   \begin{tabular}{c|l| rr|rr}
 \hline
~ & \multicolumn{1}{c|}{$N$} & \multicolumn{1}{c}{$L^\infty p$} & \multicolumn{1}{c|}{Rate} & \multicolumn{1}{c}{$L^\infty  p(t)$} & \multicolumn{1}{c}{ Rate}\\
 \hline
 \hline
 \parbox[t]{4.5mm}{\multirow{5}{*}{\rotatebox[]{90}{\parbox[c]{2.cm}{\footnotesize Adams--Moulton(4)}}}}&
 40 & 0.0220741 & 2.1699 &  7.60885e-07 & 7.26615 \\
 &80 &  0.00518945 & 2.0887 &  6.02869e-09 & 6.97969 \\
 &160 &    0.00125739 & 2.04515 &  6.24648e-11 & 6.59266 \\
 &320  &  0.000309428 & 2.02275 &  9.69648e-13 & 6.00944 \\
 &640 &  7.67471e-05 & 2.01142 &  1.69123e-14 & 5.84132 \\
 \hline
 \end{tabular}
 \end{table}

%\clearpage

\subsection{Results on the discretization of the full optimality system} 
We consider the discretization of the full optimality system \eqref{OCP} and
equations \eqref{osOCP}, respectively. Note that the example proposed in \cite{Hager00}
and also investigated in \cite{HertyPareschiSteffensen2013} is not suitable to highlight
the difference between the approaches in Figure \ref{overview} since $f_y=cst.$ Therefore,
we propose the following problem: 
\begin{equation}
\label{numOCP} 
\begin{aligned}
&\min\limits_{y,u} \frac{1}2 \left( y(T)- \frac{1}{1-T} \right)^2 + \frac{\alpha}2 \int_0^1 u^2 ds,\\
&\mbox{ subject to }  y'=y^2+u,\qquad y(0)=1,
\end{aligned}
\end{equation}
where we chose $\alpha>0$ as regularization parameter, and we remark that the exact solution for $u\equiv 0$ is given by 
$$y(t)=\frac{1}{1-t}.$$
The adjoint equations \eqref{osOCPb} and optimality conditions \eqref{osOCPc} are given by 
\begin{equation*}
\begin{aligned}
p' = 2 y p,\quad
p(T) = y(T)-\frac{1}{1-T}, \quad
p + \alpha u = 0.
\end{aligned}
\end{equation*}
Clearly, for $u=0$ we obtain $p\equiv 0.$ In order to avoid loss of accuracy due to 
inexact initialization we initialize the forward problem \eqref{osOCPa} using the 
exact solution  at time $t\leq 0$ and the adjoint equation according to the conditions 
\eqref{discOCPc}. We show the convergence results for the adjoint state $p$ 
as well as the state $y$ for different BDF methods in 
Table \ref{tab3}.

 \begin{table}[tb]
 \centering
 \caption{ BDF(4): Number of discretization points in time $N$, error in $L^\infty(0,T)$  for the approach optimize--then--discretize (Lemma \ref{lem2}) is shown in $L^\infty p$ with 
corresponding rate (Rate).  Also, shown is the $L^\infty$ error in the state $y$ in the second
column as well as its rate (Rate). We report from top to bottom different schemes: BDF(3), BDF(4), BDF(6).}  \label{tab3}
\vspace{+0.25cm}
 \begin{tabular}{c|l| rr|rr}
 \hline
~ & \multicolumn{1}{c|}{$N$} & \multicolumn{1}{c}{$L^\infty y$} & \multicolumn{1}{c|}{Rate} & \multicolumn{1}{c}{$L^\infty  p$} & \multicolumn{1}{c}{ Rate}\\
 \hline
 \hline
 \parbox[t]{4.5mm}{\multirow{6}{*}{\rotatebox[]{90}{\parbox[c]{1.6cm}{\footnotesize BDF(3)}}}} &40 & 0.0720175 & 2.94884 &  3.47941 & 3.44822 \\
 &80 & 0.0107919 & 2.73839 &  0.498712 & 2.80257  \\
 &160 & 0.00153707 & 2.8117 &  0.0705343 & 2.82181  \\
 &320 & 0.000207256 & 2.8907 &  0.00950082 & 2.8922  \\
 &640 & 2.6974e-05 & 2.94177 &  0.00123634 & 2.94198  \\
 &1280 & 3.44239e-06 & 2.97009 &  0.000157777 & 2.97011 \\
  \hline
% \end{tabular}
%\\
%~
%\\
%  \begin{tabular}{c|l| rr|rr}
 \hline
~ & \multicolumn{1}{c|}{$N$} & \multicolumn{1}{c}{$L^\infty y$} & \multicolumn{1}{c|}{Rate} & \multicolumn{1}{c}{$L^\infty  p$} & \multicolumn{1}{c}{ Rate}\\
 \hline
 %\hline
  \parbox[t]{4.5mm}{\multirow{6}{*}{\rotatebox[]{90}{\parbox[c]{2.cm}{\footnotesize BDF(4)}}}}& 40 & 0.0237103 & 3.32788 &  1.25952 & 3.56845  \\
 &80 & 0.00224529 & 3.40054 &  0.117177 & 3.42611  \\
 &160 & 0.000182662 & 3.61966 &  0.00951526 & 3.6223  \\
 &320 & 1.32309e-05 & 3.78719 &  0.000689121 & 3.78741  \\
 &640 & 8.93826e-07 & 3.88778 &  4.65535e-05 & 3.8878  \\
 &1280 & 5.81525e-08 & 3.94208 &  3.02878e-06 & 3.94208  \\
 \hline
% \end{tabular}
%\\
%~
%\\
%   \begin{tabular}{c|l| rr|rr}
 \hline
~ & \multicolumn{1}{c|}{$N$} & \multicolumn{1}{c}{$L^\infty y$} & \multicolumn{1}{c|}{Rate} & \multicolumn{1}{c}{$L^\infty  p$} & \multicolumn{1}{c}{ Rate}\\
 \hline
 \hline
 \parbox[t]{4.5mm}{\multirow{6}{*}{\rotatebox[]{90}{\parbox[c]{2.cm}{\footnotesize BDF(6)}}}}
& 40 & 0.00451569 & 4.10057 &  0.27787 & 4.19135  \\
 &80 & 0.000188028 & 4.58593 &  0.0115192 & 4.59229  \\
 &160 & 5.42671e-06 & 5.11473 &  0.000332388 & 5.11503  \\
 &320 & 1.20044e-07 & 5.49844 &  7.35271e-06 & 5.49845  \\
 &640 & 2.2528e-09 & 5.7357 &  1.37984e-07 & 5.7357  \\
 &1280 & 2.40865e-11 & 6.54735 &  1.4753e-09 & 6.54735  \\
 \hline
 \end{tabular}
 \end{table}

%\clearpage

\subsection{ BDF discretization for the relaxation system and adjoint}

In this section we consider the discretized relaxation system \eqref{ODE} being the forward problem as well as the corresponding discretized adjoint equation given by equation \eqref{implement adjoint}. 

%In Table \ref{tab4} we show convergence rates for the 
%macroscopic quantity $u(t,x)=\frac{1}N\sum_j f_j$ computed using 
%BDF discretizations of equation \eqref{ODE}. We consider the case
%$N=2$ and $v_1=-v_2=a=2.1$ and the two different test cases of
%pure advection, $f(u)=u$, and Burger's equation $f(u)=\frac{u^2}2.$
%The initial data is $u_0(x)=\exp( - (x-3)^2)$ 
%and terminal  time is $T=1$ on a domain $x \in [0,6]$ with periodic boundary conditions for both cases. Numerical grids with $N_x=20$ to $N_x=640$ points are considered. 
%The value of $\epsilon$ is kept fixed at $\epsilon=10^{-2}.$ The temporal grid
%is chosen according to the CFL condition $\Delta t = \Delta x$. 

\paragraph{Forward system.} We study numerically the evolution of the macroscopic quantity $u(t,x)=\frac{1}N\sum_j f_j$ computed using 
BDF discretization of equation \eqref{ODE}. We consider the case
$N=2$ and $v_1=-v_2=a=2.1$ and two different test cases of
pure advection, $F(u)=u$, and Burger's equation $F(u)=\frac{u^2}2.$
The initial data is $u_0(x)=\exp( - (x-3)^2)$  and terminal  time is $T=1$ on a domain $x \in [0,6]$ with periodic boundary conditions for both cases. We considered $N_x=640$ grid points for the space discretization, and the temporal grid
is chosen according to the CFL condition, such that $\Delta t = \Delta x$, the value of $\epsilon$ is kept fixed at $\epsilon=10^{-2}$.

\begin{figure}[tb]
\centering
\ $F(u)=u$ \hspace{+4.5cm} $F(u)=u^2/2$ \
\\
~
\\
\includegraphics[scale = 0.35]{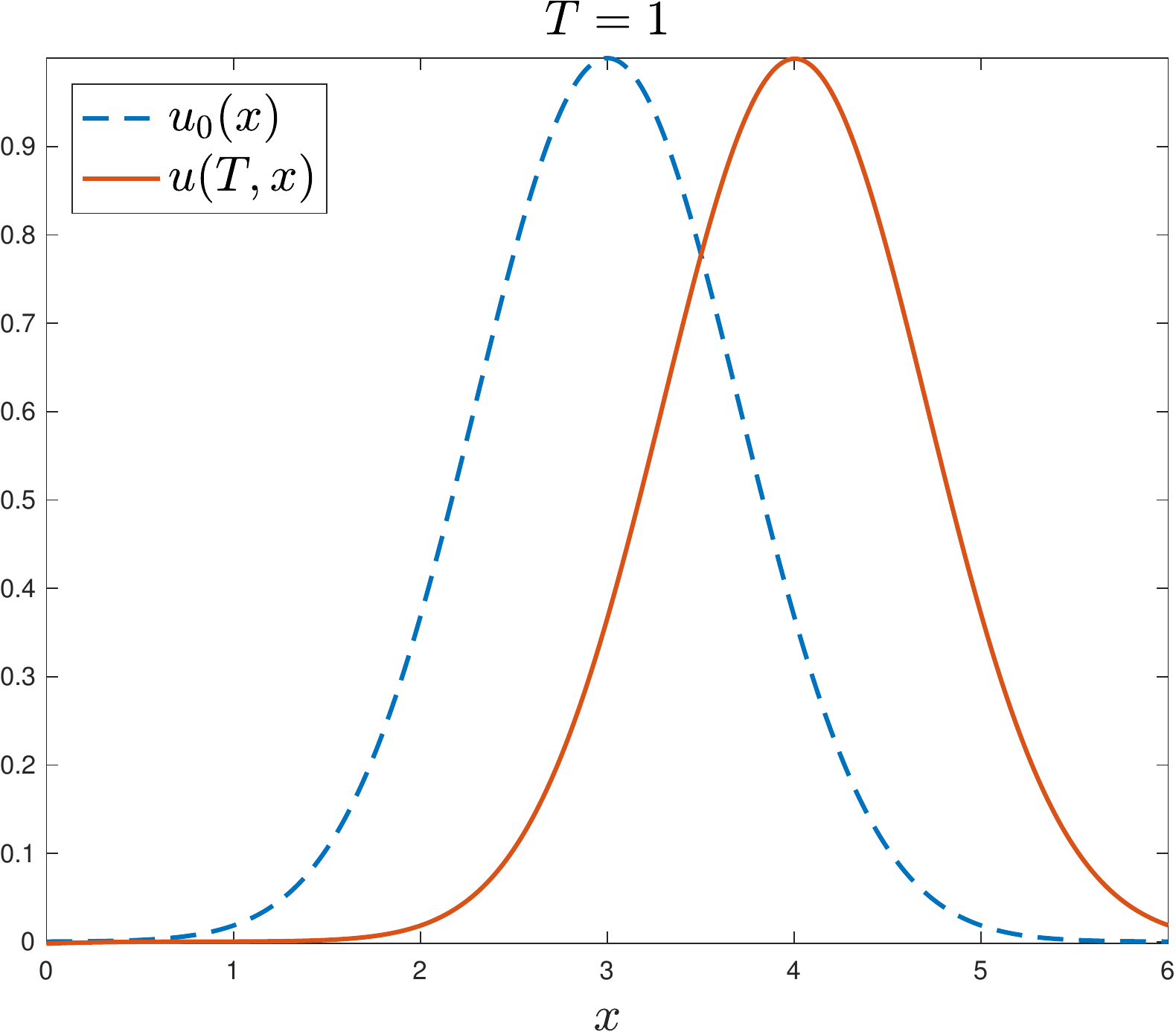}
\qquad\qquad
\includegraphics[scale = 0.35]{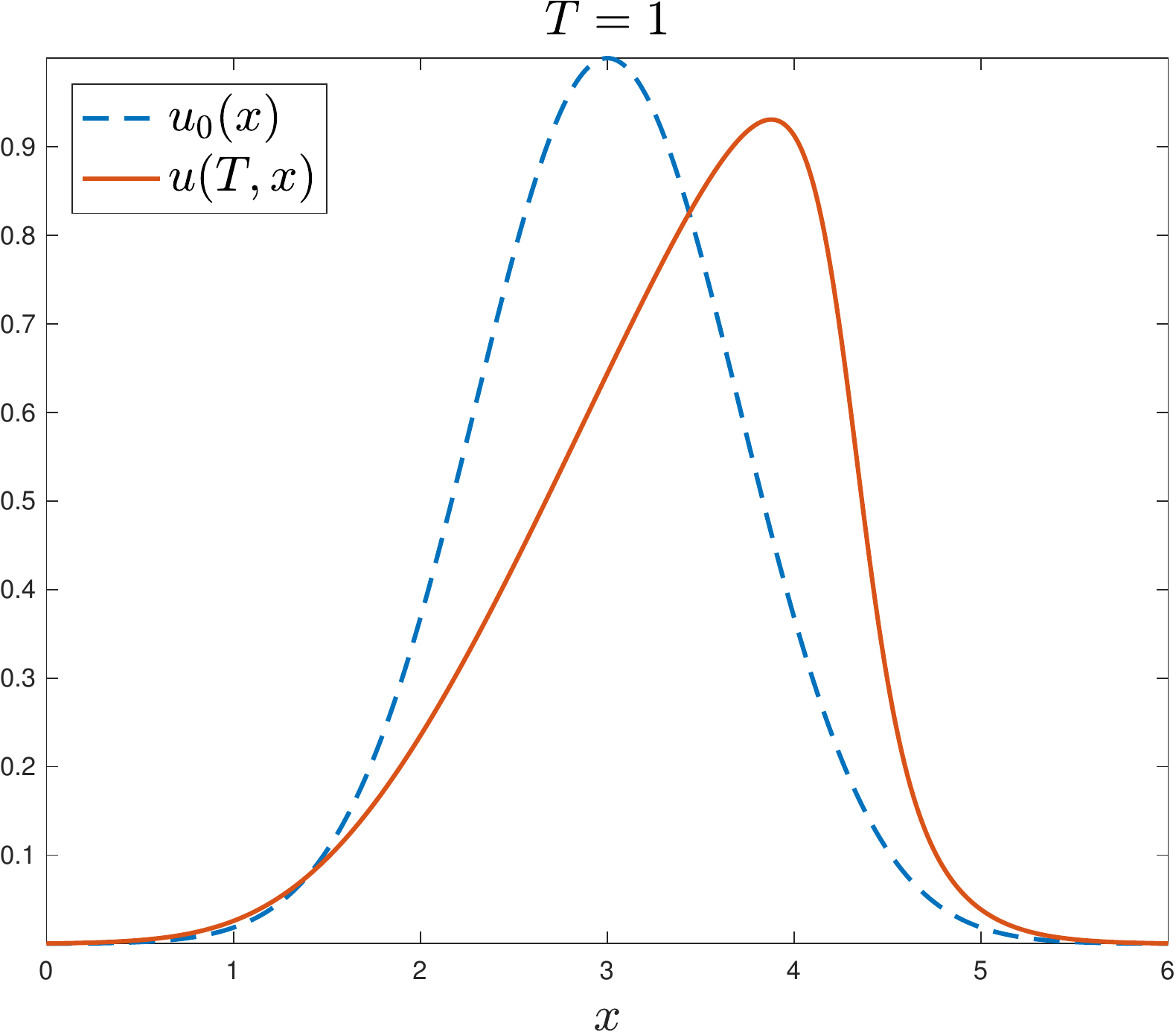}
\\
\includegraphics[scale = 0.35]{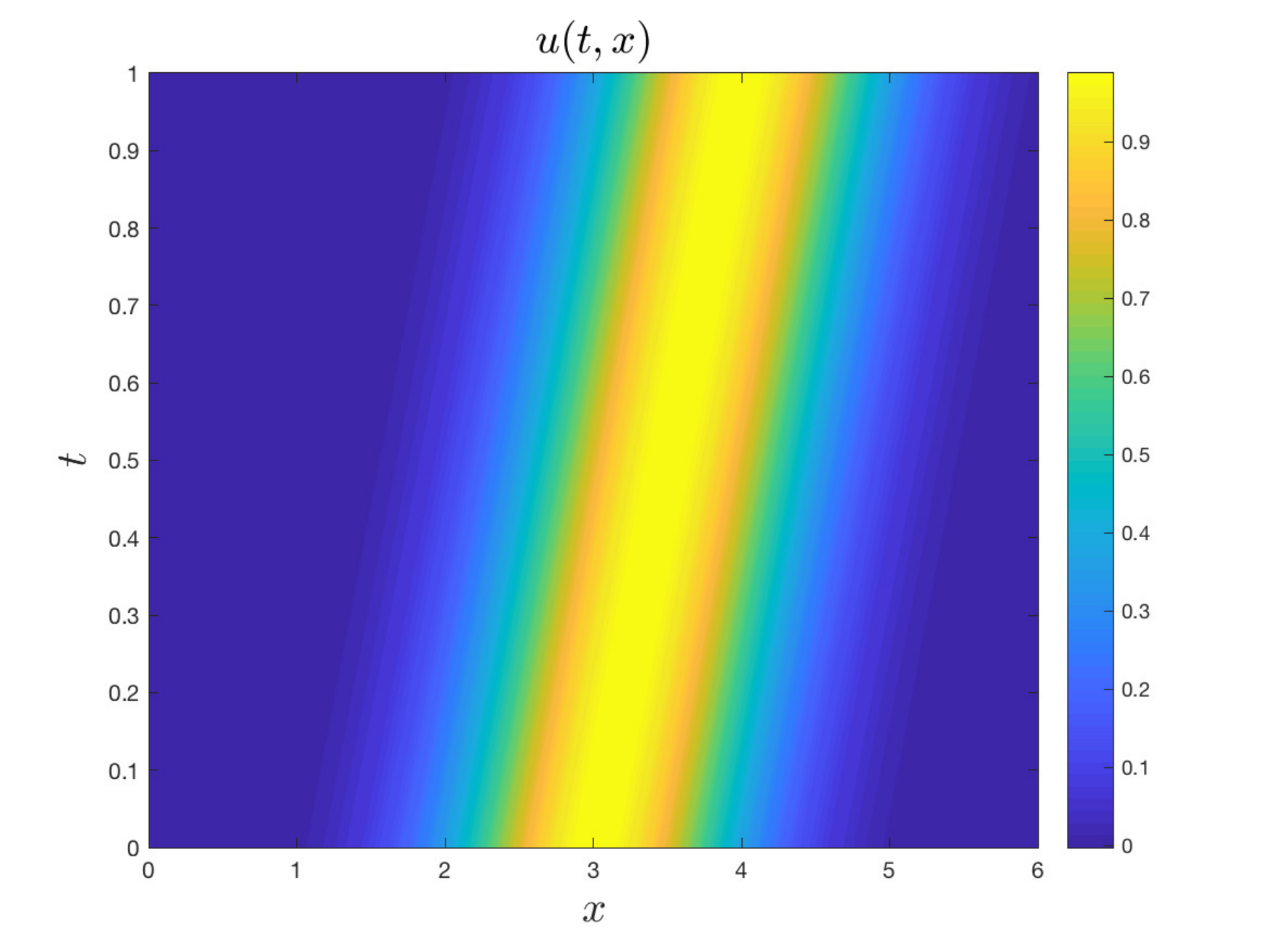}
\,
\includegraphics[scale = 0.35]{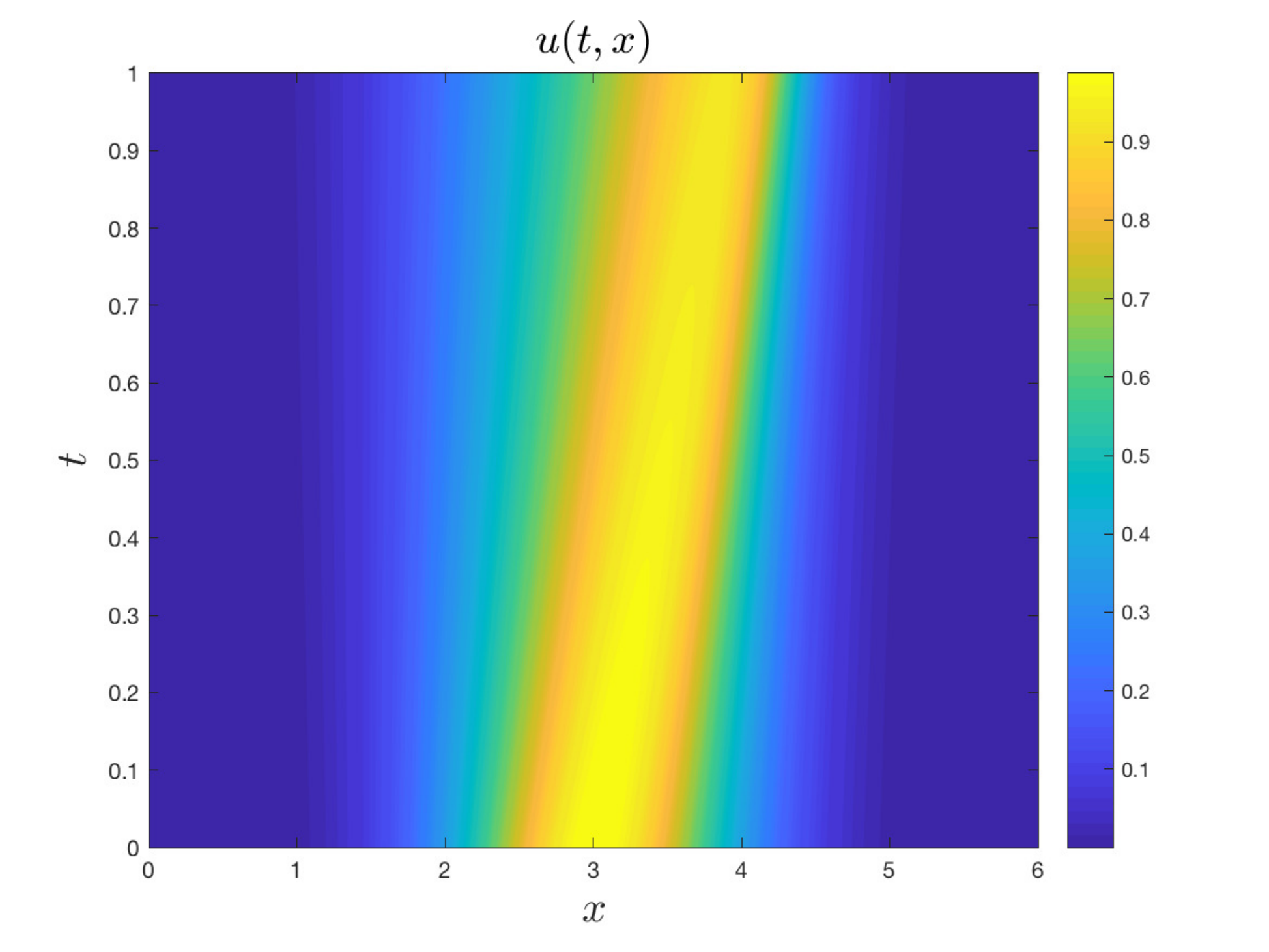}
\caption{BDF integration of the system of ODEs \eqref{ODE} used as BGK approximation to the conservation law \eqref{cons law}. Two velocities are considered, $N=2$.  Left-hand side column corresponds to pure transport situation $F(u) = u$, whether the right-hand side column depict the solution of the Burger flux function, $F(u) = u^2/2$. Top row  show initial data $u_0(x)$ as well as the numerical result at terminal time $T =1$, bottom row shows the density $u(x,t)$ in the space-time frame $[0,6]\times[0,1]$. Each test has been produced using a BDF(3) scheme with $N_x  =640$ space points and $\dt = 4.47127\times 10^{-3}$, and with fixed relaxation parameter $\epsilon = 0.01$.}\label{Fig001}
\end{figure}

We present the numerically solutions in Figure \ref{Fig001} for the linear and non-linear transport case.
Here, higher-order successfully reduces the numerical diffusion and yields qualitatively better results. 
\clearpage
We do not present convergence tables for the forward equation since  equation \eqref{ODE} requires to evaluate the local equilibrium at 
gridpoints $x+v_j t$ that are in general not aligned with the numerical 
grid. Therefore, an interpolation is required. Hence, the temporal and spatial resolution 
are not independent and the observed convergence is limited to the interpolation 
of the solution. 

\begin{figure}[tb]
\centering
\centering
\includegraphics[scale = 0.245]{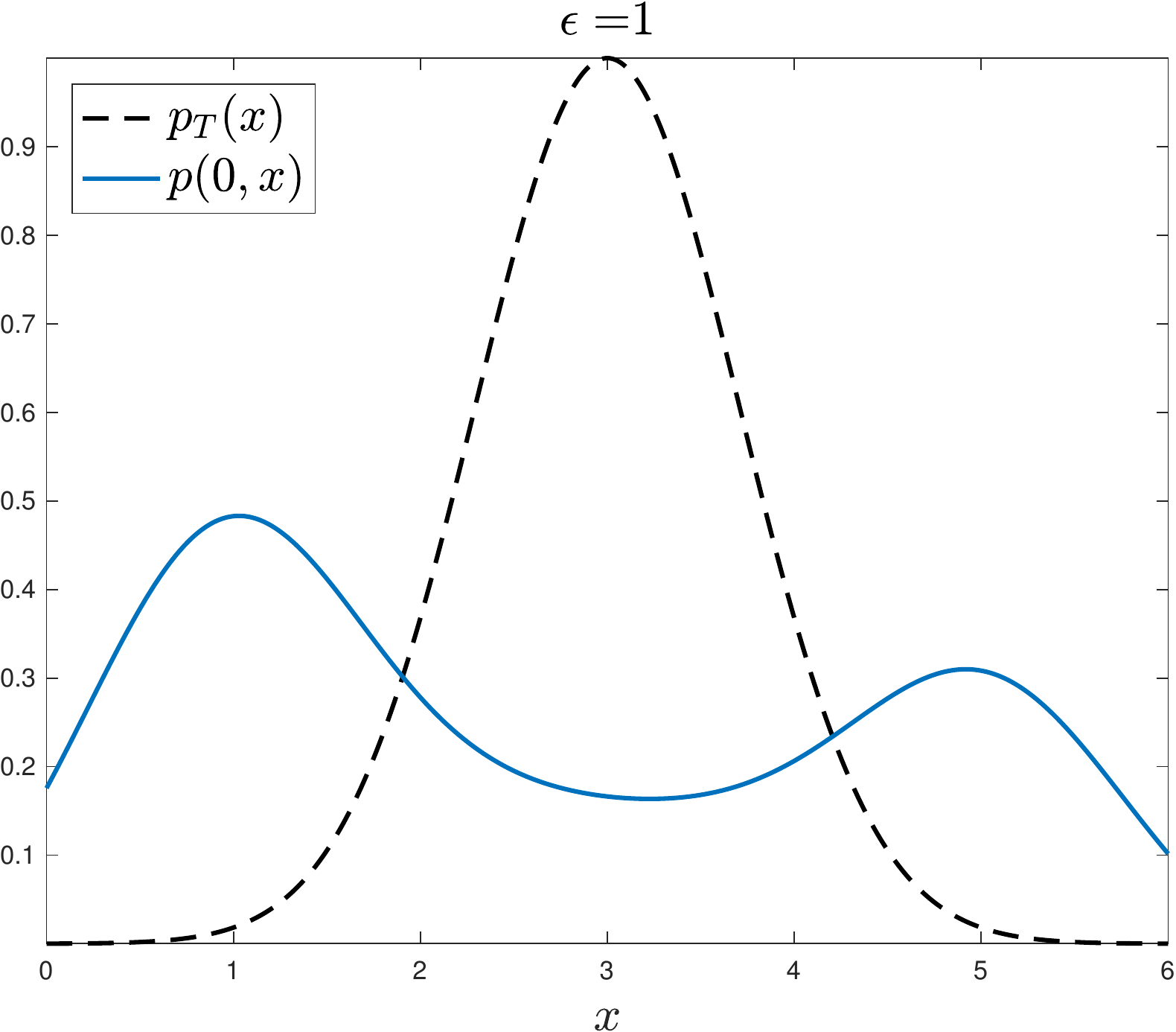}
\quad\qquad
\includegraphics[scale = 0.245]{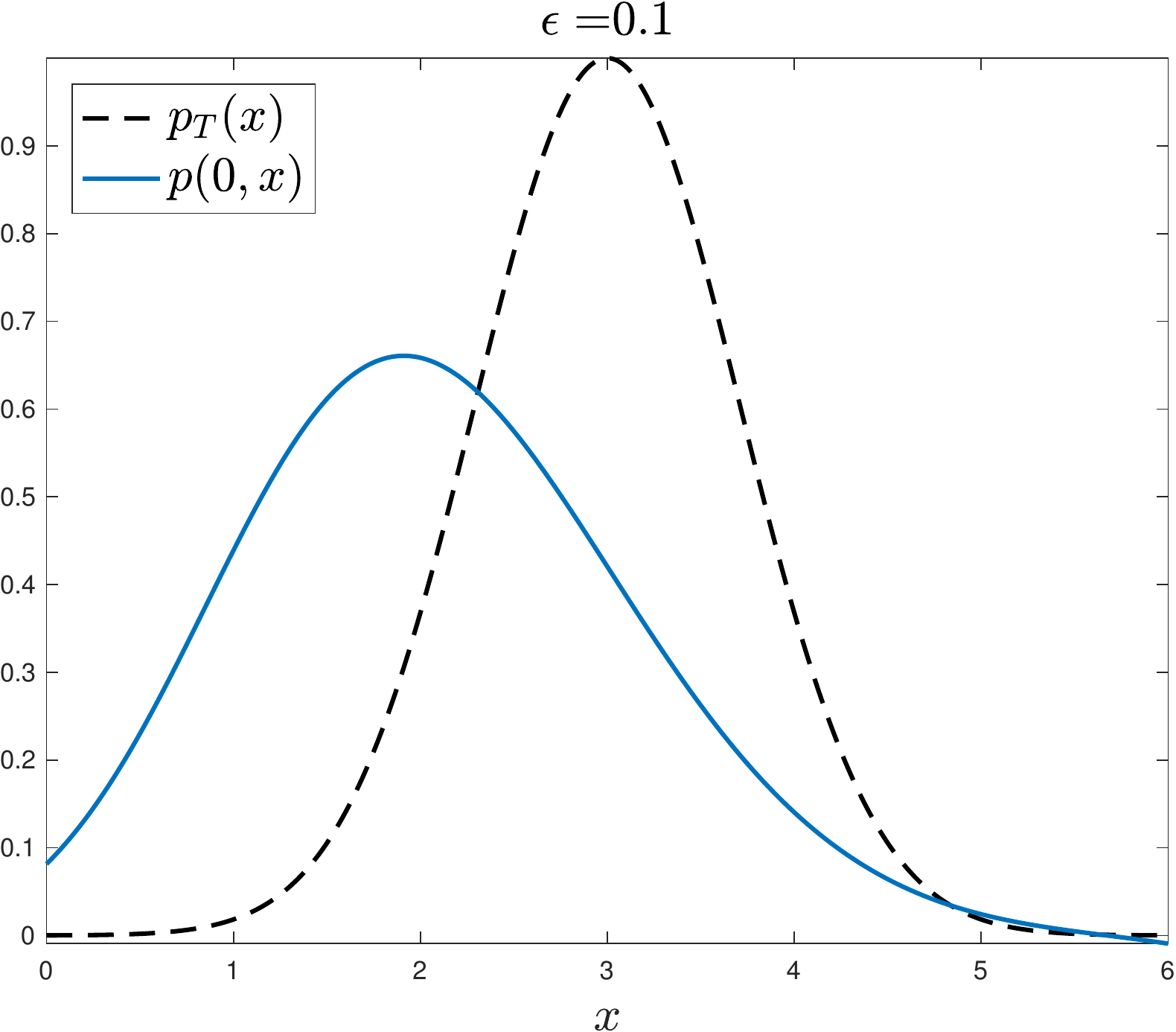}
\quad\qquad
\includegraphics[scale = 0.245]{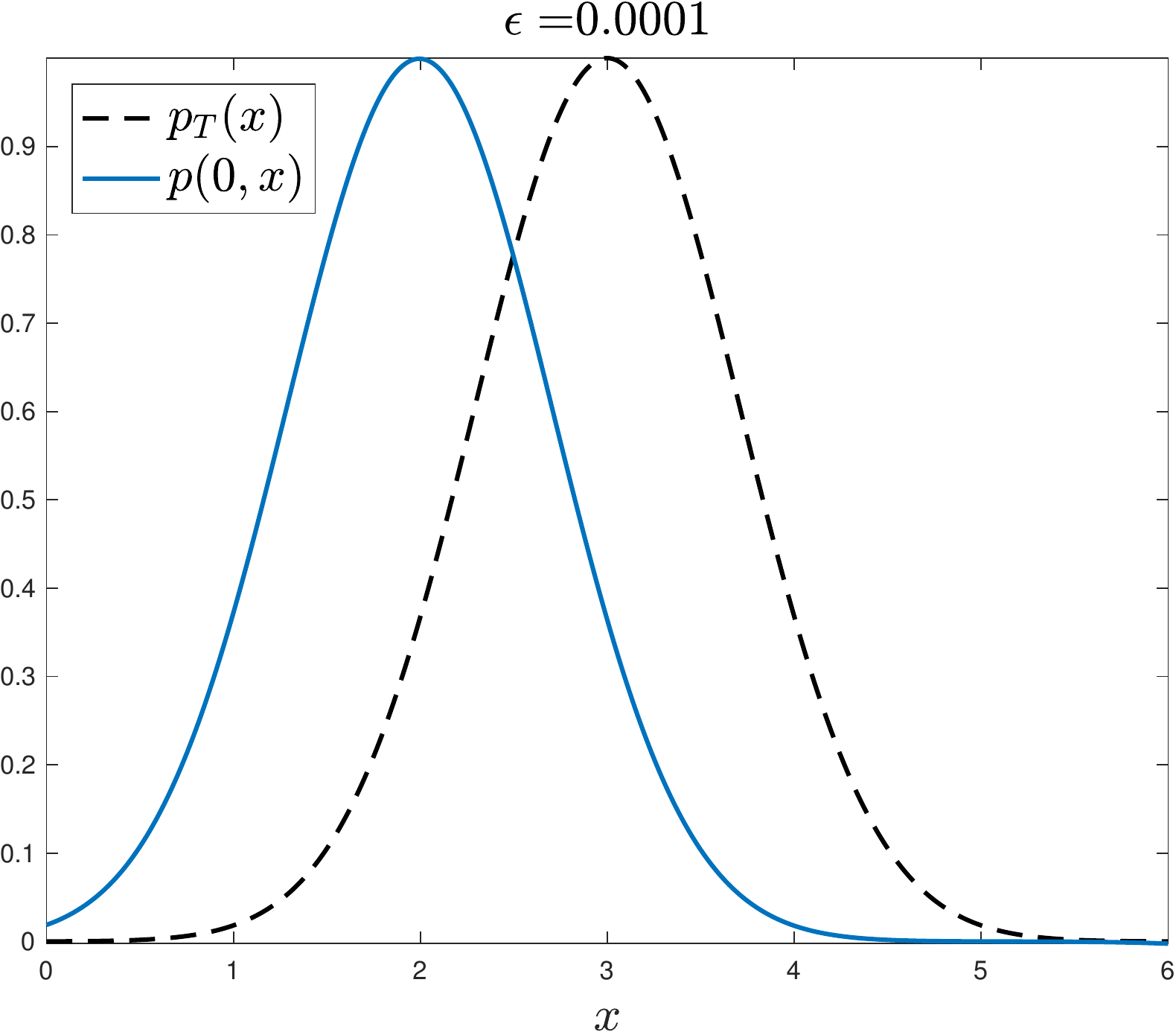}
\\
\includegraphics[scale = 0.245]{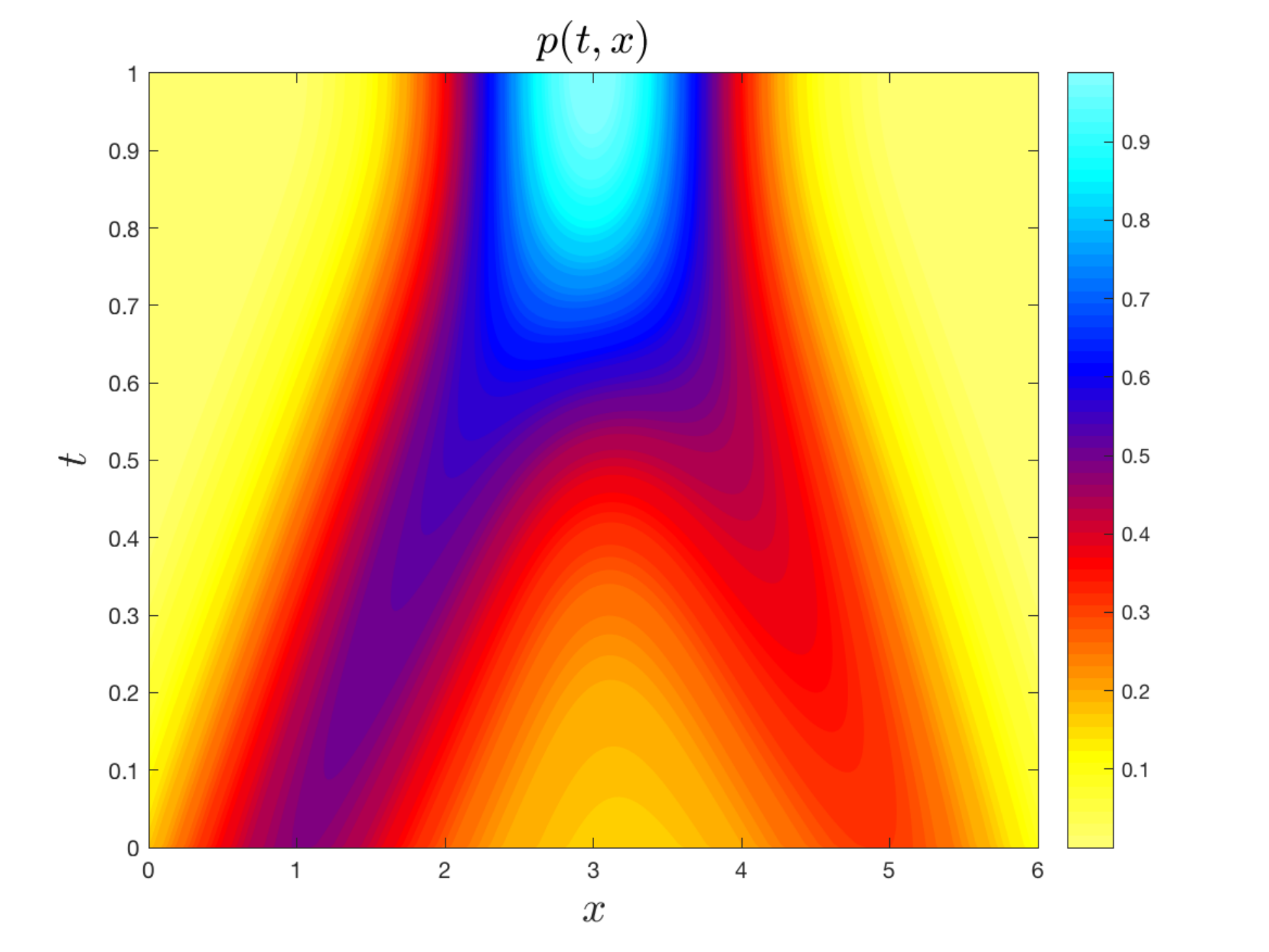}
\includegraphics[scale = 0.245]{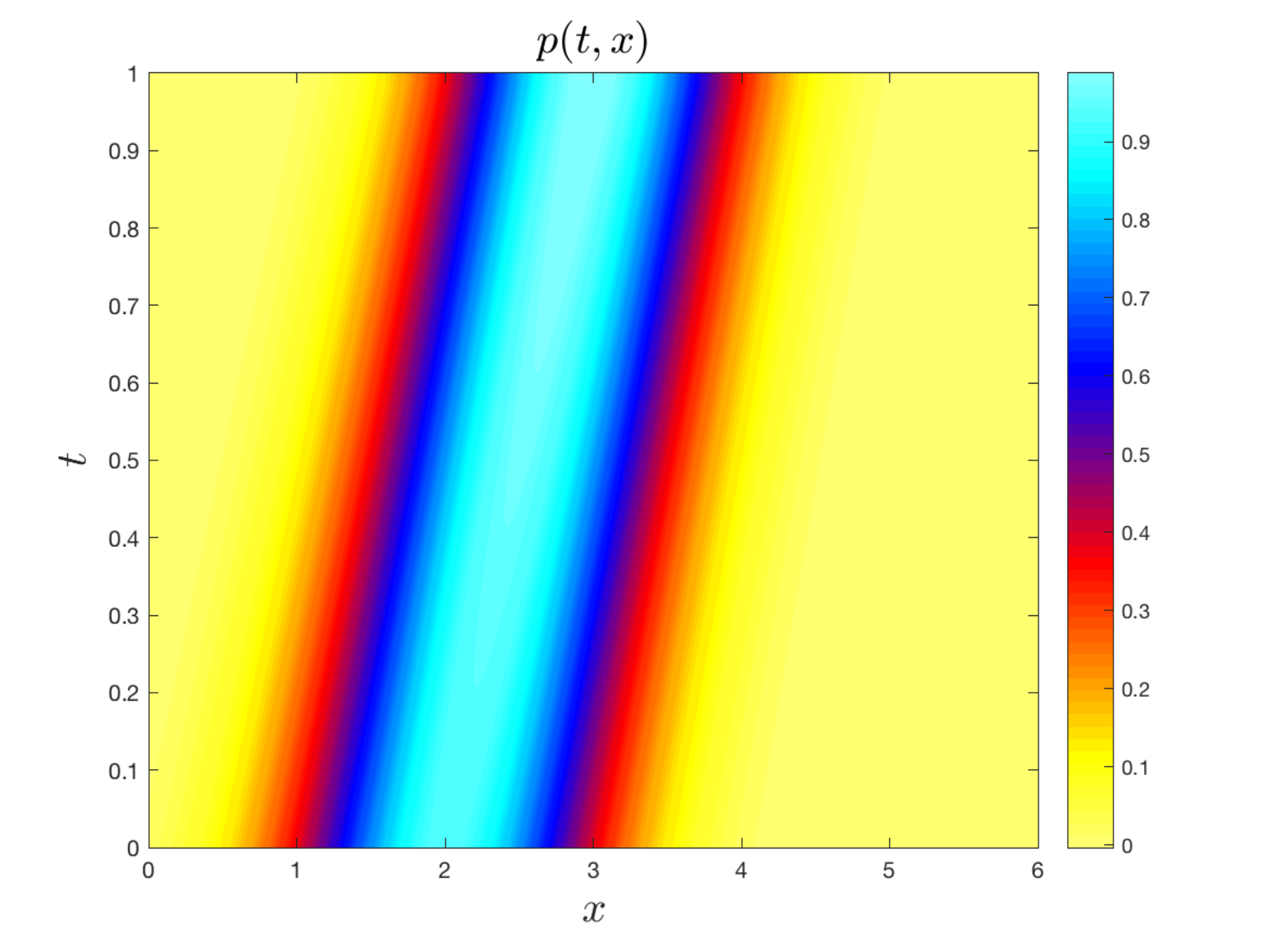}
\includegraphics[scale = 0.245]{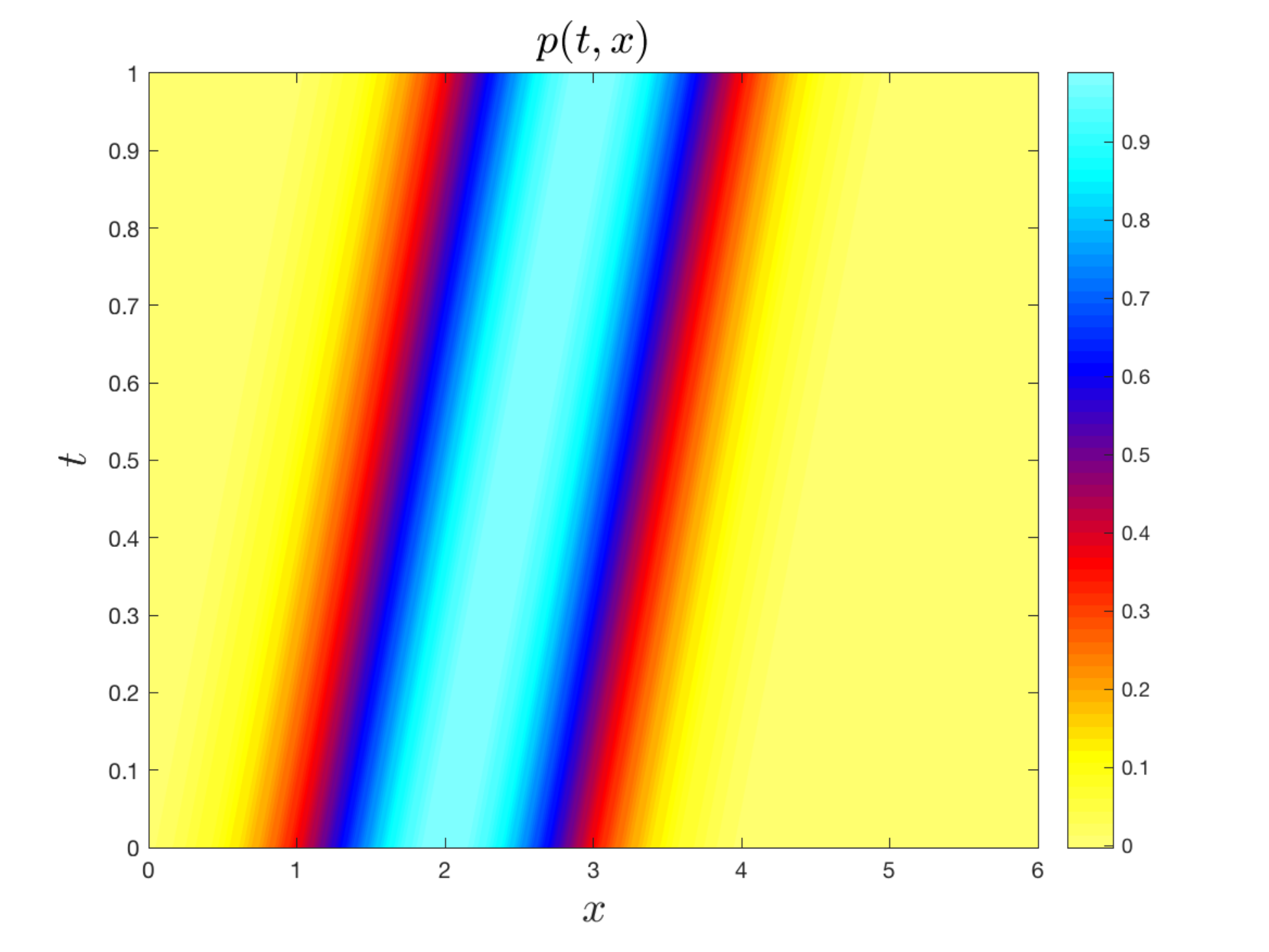}
\caption{The BDF(2) integration of the system of \eqref{implement adjoint} has been implemented for the linear transport $F(u)=u.$ Two velocities are considered,  with $N_x  =640$ space points and $\dt = 4.47127\times 10^{-3}$. From left to right we show different values of $\epsilon$, with $\epsilon\in{\{1,0.1,0.0001\}}$. In the top row the terminal data $p_T(x)$, $T=1$, as well as the numerical result at initial time $p(0,x)$ are reported for the different values of $\epsilon$. Bottom row depict the density $p(t,x) = \lambda^1(t,x)+\lambda^2(t,x)$ for the different value of $\epsilon$. Note that for small $\epsilon$ the pure transport equation is obtained.}\label{Fig002}
\end{figure}

\paragraph{Adjoint system.}
A similar behavior is observed for the discretization 
of the adjoint equation \eqref{BGK adjoint}. In order to 
illustrate the results we only show the BDF(2) method applied 
to \eqref{implement adjoint} in the case of $F(u)=u.$ We use
the same parameters as above for the forward system, but now the data $u_0(x)$ 
is prescribed at terminal time $T=1$, in the following way 
$\lambda^j(T,x)=u_0(x)/N$, with $N = 2$. 
Then, the adjoint variables are evolved according to the derived scheme \eqref{implement adjoint}. 
For illustration purposes the solutions $p(t,x) = \lambda^1(t,x)+\lambda^2(t,x)$ are reported for different values of the scaling term $\epsilon$ in Figure \ref{Fig002}, in the top row  we represent  the adjoint equation at time zero jointly with  the terminal conditions $p_T(x)$, in the bottom row the density $p(t,x)$ in the domain $[0,6]\times[1,0]$.
Compared with the Figure \ref{Fig002} we observe that the profile moves over time in the opposite direction, when $\epsilon$ is small enough.
This is precisely as expected by the limiting equation $-p_t - F'(u) p_x =0$, where $p = \sum_j \lambda^j.$ 
%\begin{figure}[h]
%\centering
%\centering
%\includegraphics[scale = 0.35]{l0lT_bwdlinear2.pdf}
%\qquad\qquad
%\includegraphics[scale = 0.35]{l0lT_bwdlinear5.pdf}
%\\
%\includegraphics[scale = 0.35]{l0T_bwdlinear2.pdf}
%\,
%\includegraphics[scale = 0.35]{l0T_bwdlinear5.pdf}
%\caption{BDF integration of the system of \eqref{implement adjoint}. Two velocities are considered and various numerical grids displayed in different colors. We show the terminal data as well as the numerical result at initial time. The finest numerical discretization as well as the initial time is displayed in the title of the corresponding subfigure.  The BDF(2) scheme has been applied
%for $f(u)=u.$}\label{Fig002}
%\end{figure}

Finally, we study the dependence of the  adjoint equation on the parameter $\epsilon.$ Note that for equation \eqref{ODE}  a similar study has been performed in \cite{GroppiRussoStracquadanio}. For each fixed value of $\epsilon$  we compute the average converge rate on the numerical grid given above. 
We also record the minimal error as well as the minimal used time step.  The study is done for the BDF(2) scheme and the results are reported
in Table \ref{tab04}.
\begin{table}[tb]
\center
 \caption{BDF integration of the system of ODEs \eqref{BGK adjoint}. Two velocities are considered and the $L^2$ error of $p= \lambda^1+\lambda^2$ at initial time is reported.  Various values of $\epsilon$ are considered. 
 The mean convergence rate on given grid is reported as well as the 
 finest temporal grid size considered. }\label{tab4}
 \vspace{+0.25cm}
\begin{tabular}{l|l|l|l}
\hline
$\epsilon$ & $\dt=\Delta x$ & $L^2p$ & Rate \\
 \hline
 \hline
1 & 0.00447127 & 8.75403e-06 & 2.74404 \\
 1.000000e-01 & 0.00447127 & 6.63095e-06 & 2.68572 \\
 1.000000e-02 & 0.00447127 & 1.77628e-05 & 2.62852 \\
 1.000000e-03 & 0.00447127 & 2.08908e-05 & 2.6135 \\
 1.000000e-04 & 0.00447127 & 2.12566e-05 & 2.61174 \\
 \hline
 \end{tabular}
 \label{tab04}
 \end{table}
%\newpage
\subsection{Optimal control of  hyperbolic balance laws}
We finally show the quality of our approach by two applications to the optimal control of hyperbolic balance laws. For further references, and example about optimal control problems governed by conservation laws we refer to \cite{CHK,HKK15}.

\subsubsection{Jin-Xin relaxation system}
We consider the Jin-Xin relaxation model, \cite{JX} which results in the two velocities model \eqref{bgk}, as follows
\begin{equation}\label{eq:realxBurger}
\begin{split}
 f^{(1)}_t + a f^{(1)}_x &= \frac{1}\epsilon \left( E_1(u) -f^{(1)} \right),\qquad f^{(1)}(x,0) = f^{(1)}(x) 
 \\ 
 f^{(2)}_t - a f^{(2)}_x &= \frac{1}\epsilon \left( E_2(u) - f^{(2)} \right),\qquad f^{(2)}(x,0) = f^{(2)}(x)
\end{split}
\end{equation}
where  the equilibrium states $E_{1}$ and $E_2$ are given by  
$$ E_{1}(u) = \frac{1}{2a} \left( au + F(u) \right), \quad E_2(u) = \frac{1}{2a} \left(au- F(u) \right),$$
with total density $u = f^{(1)} + f^{(2)} $ and velocity in the limit $\epsilon\to0$ such that $v=a(f^{(1)}-f^{(2)})= F(u)$. In particular we choose as an example the flux $F(u) = u^2/2$ associated to the inviscid Burger equation.
Thus we have that the characteristic speed $a$ has to satisfy the condition $ a \geq \max_{ x \in  \R } | u_0(x)|.$ 
We consider the following optimal control problem, firstly proposed in \cite{HKK15}, where here we seek for minimizers $f^{(j)}_0(x)$, with $j=1,2$ of
\begin{align}\label{optJu}
J(u(\cdot,T),u_d(\cdot)) = \frac{1}{2}\int_{\Omega}|u(x,T)-u_d(x)|^2\, dx
\end{align}
%subject to the evolution of the invisicid Burger equation defined as follows
%\begin{equation}\label{burger}
%\begin{aligned}
%&u_t +\left(\frac{u^2}{2}\right)_x = 0,\qquad u(0,x) = u_0,\qquad (t,x)\in[0,T]\times\Omega,
%\end{aligned}
%\end{equation}
Hence we  fix the specific domain $\Omega = [-3,3]$ with $T=3$ and we want to prescribe the final discontinuous data $u_d(x)$, as final data at time $T$ of the Burger equation with initial data defined as follows
\begin{equation}\label{u_d}
\begin{aligned}
u_d(0,x) = 
\begin{cases}
1.5+x\qquad \textrm{if}\quad -1.5\leq x\leq-0.5,\\
0\qquad  \textrm{otherwise.}
\end{cases}
\end{aligned}
\end{equation}
%Then an exact solution of optimal control problem  \eqref{optJu}--\eqref{eq:realxBurger} can be derived explicitly, see \cite{HKK15}, and it is defined as follows
%\begin{equation}\label{uopt}
%\begin{aligned}
%\hat{u}_0(x) = 
%\begin{cases}
%1&\quad \textrm{if}\quad -1/2\leq x<0\\
%1-x&\quad \textrm{if}\quad 0\leq x\leq1\\
%0&\quad  \textrm{otherwise}
%\end{cases}
%\qquad
%{u}(T,x) = 
%\begin{cases}
%x+\frac{1}{2}&\quad \textrm{if}\quad -1/2\leq x<1/2\\
%1&\quad \textrm{if}\quad 1/2\leq x\leq1\\
%0&\quad  \textrm{otherwise}
%\end{cases}
%\end{aligned}
%\end{equation}
In order to solve numerically this optimal control problem we approximate it with the optimal control \eqref{OPTx}--\eqref{BGKforward}, choosing $N=2$ velocities and $N_x  =120$ space points, time step $\dt = 0.05$ and relaxation parameter $\epsilon =10^{-2}$. In order to solve the time discretization we use BDF(2) integration. The same choice of parameters is considered for the adjoint equation \eqref{BGK adjoint}, which is solved backward in time using a two velocities approximation of the terminal condition $\lambda^{(1)}(T,x)+\lambda^{(2)}(T,x)=p(T,x) = u(x,T)-u_d(x)$. 
Thus, we solve recursively the forward approximated system \eqref{BGKforward} and the backward system \eqref{BGK adjoint}, 
using as starting point the step function $u^{(0)}_0(x) = 0.5\chi_{[-1.5,-0.5]}(x)$, and introducing a filter $\mathcal{F}(\cdot)$ to reduce the total variation of the initial data $u^{(k)}_0$ following the approach in \cite{HKK15}. Then we update the initial data $u_0^{(0)}(x)$ using a steepest descend method as follows
\[
u^{(k+1)} _0=u_0^{(k)}-\sigma_k p^{(k)}, \qquad m\geq 0,
\]
with $\sigma_k$ updated with Barzilai-Browein step method \cite{BB}.

\begin{figure}[tb]
\centering
\centering
\includegraphics[scale = 0.35]{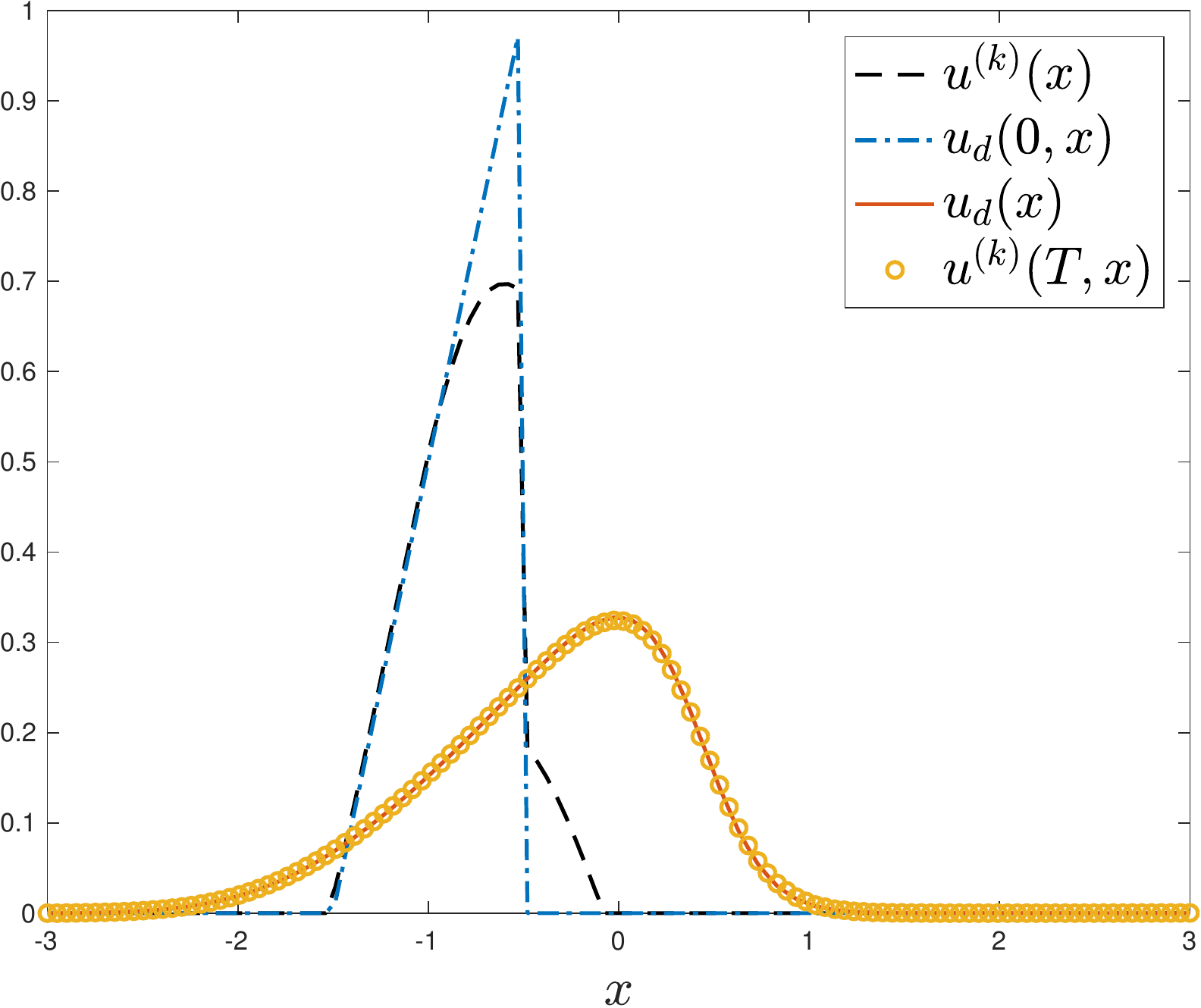}
\quad\qquad
\includegraphics[scale = 0.35]{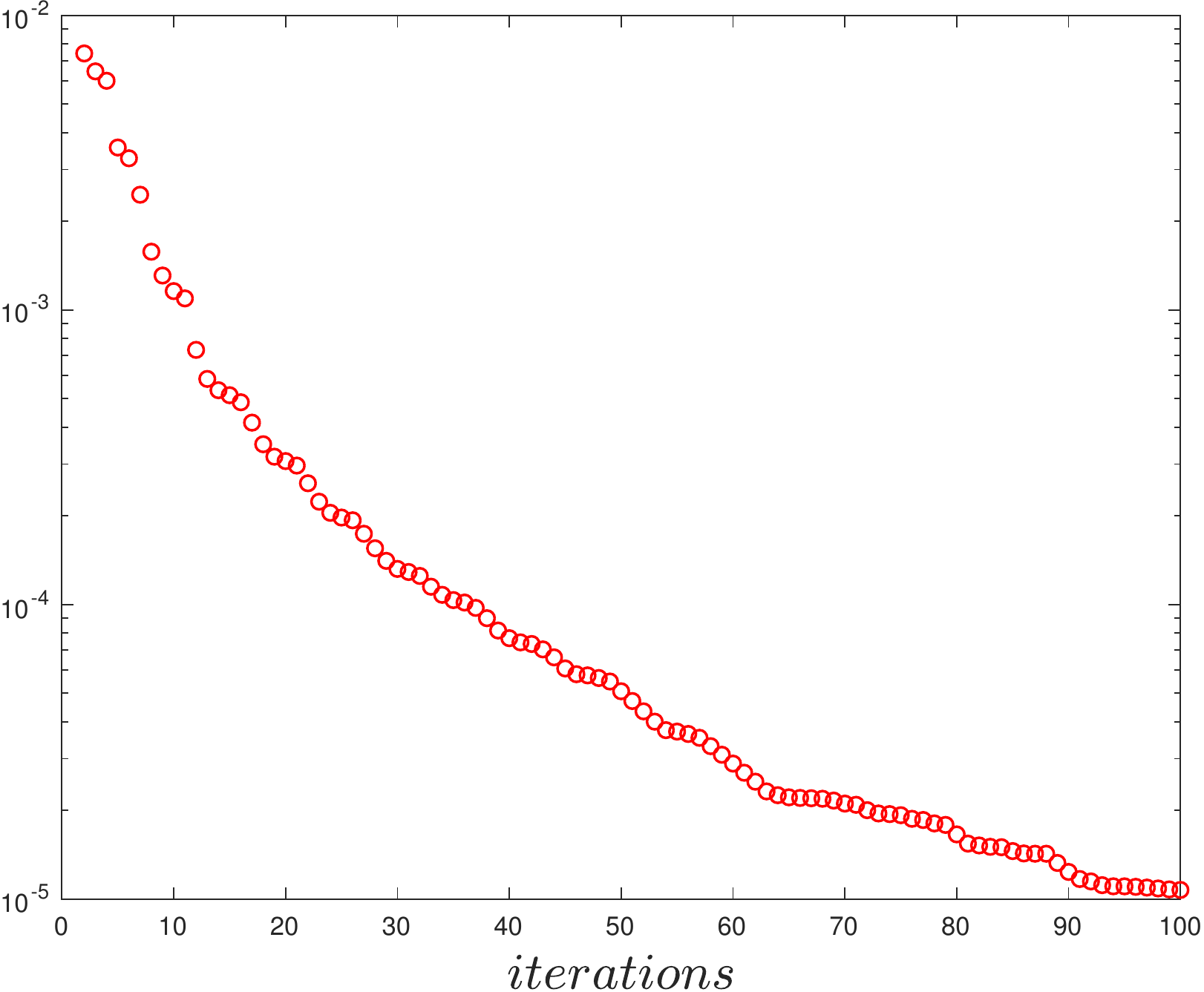}
\caption{{\em Jin-Xin relaxation system}. On the left-hand side we report the control at iteration 100, $u^{(100)}_0$ (\texttt{--}), compared with the initial data $u_d(0,x)$ (\texttt{-.}) used to compute the desired final data $u_d(x)$ ($\texttt{-}$). The terminal solution $u(T,x)$, (\texttt{o}), is reported and it is computed solving system \eqref{eq:realxBurger}  using BDF(2) integration, with $N_x  =120$ space points, time step $\dt = 0.05$, and relaxation parameter $\epsilon =10^{-2}$. On the right-hand side we show the decrease of the functional $J(u^{(k)})$.}\label{Fig004}
\end{figure}

We report in Figure \ref{Fig004} the final result after $k = 30$ iterations of the optimization process, on the left-hand side plot we depict the initial data $u^{(100)}_0$ with the terminal data $u(T,x)$ as well as the desired data $u_d(x)$. On the right hand side we depict the decrease of $J(u^{(k)})$given by the optimization procedure.
%Notice that the shape of the $u^{(30)}_0$ is a smooth approximation of the exact solution \eqref{uopt}.  

%In Figure \ref{Fig005} we report the decrement of the functional \eqref{OPTx} accordingly to each value of $u^{(k)} _0$ with $k$, and stopping the optimization after $k=30$ iterations.

\subsubsection{Broadwell model}
We consider the one-dimensional Broadwell model, \cite{Broadwell}, which describe the evolution of densities $f^{(1)},f^{(2)},f^{(3)}$ relative to the velocities $c,-c,0$, with $c>0$, as follows
\begin{align}
\begin{aligned}\label{eq:relaxIE}
&f^{(1)}_t + c f^{(1)} _x =\frac{1}{\epsilon} \left(E_1(\rho,m) - f^{(1)}\right),\qquad f^{(1)}(x,0) = f^{(1)}_0(x),  \\
&f^{(2)}_t - c f^{(2)} _x =\frac{1}{\epsilon} \left(E_2(\rho,m) - f^{(2)}\right),
\qquad f^{(2)}(x,0) = f^{(2)}_0(x),  \\
&f^{(3)}_t  =\frac{1}{\epsilon} \left(E_3(\rho,m) - f^{(3)}\right),\quad\qquad\qquad f^{(3)}(x,0) = f^{(3)}_0(x).
\end{aligned}
\end{align}
Where the equilibrium quantities are defined as follows
\begin{align*}
E_1(\rho,m,z) &=\frac{1}{2}F(\rho,m)+\frac{m}{2c},\\E_2(\rho,m,z) &=\frac{1}{2}F(\rho,m)-\frac{m}{2c},\\E_3(\rho,m,z) &=-F(\rho,m)+\rho
\end{align*}
and the macroscopic quantities $m, \rho$, jointly with the flux $F(\rho,m)$ are such that
\begin{align}\label{eq:trasf}
 \rho = f^{(1)}+f^{(2)}+2f^{(3)},\qquad m = c(f^{(1)}-f^{(2)}),\qquad F(\rho,m) = \frac{m^2}{c^2\rho} +\rho.
 \end{align}
Indeed for $\epsilon\to0$ system \eqref{eq:relaxIE} converges to to the isentropic Euler model, \cite{HKK15}, where  $\rho, m$ represent respectively density, and momentum, 
\begin{equation}\label{eq:IE}
\begin{cases}
\rho_t + m _x = 0,\\
m_t + c^2\left(\frac{m^2}{c^2\rho} +\rho\right)_x = 0,\qquad (x,t)\in\mathbb{R}\times(0,T]
\end{cases}
\end{equation}
We aim to minimize the functional 
\begin{equation}\label{eq:JIE}
J(\rho(\cdot,T),m(\cdot,T)) = \frac{1}{2}\int_\Omega\left( |\rho(x,T)- \rho_d(x)|^2 + |m(x,T)- m_d(x)|^2  \right) \ dx
\end{equation}
with respect to the initial data $f^{j}_0(x)$ for $j=1,2,3$ taking in to account the relations \eqref{eq:trasf}.
%In order to introduce a relaxation approximation of system \eqref{eq:IE} we consider the variable $z$, and the hyperbolic relaxation system with parameter $\epsilon>0$ as follows
%\begin{equation}\label{eq:IErel}
%\begin{cases}
%\rho_t + m _x = 0,\\
%m_t +c^2 z_x = 0,\\
%z_t + m_x =\frac{1}{\epsilon}\left(F(\rho,m)-z\right),
%\end{cases}
%\end{equation}
%such that for $\epsilon\to0$ the system relaxes to the equilibrium $F(\rho,m) = m^2/(c^2\rho) +\rho$.
%
%Hence from \eqref{eq:IErel} we obtain the three velocities Broadwell model as follows
%\begin{align}
%\begin{aligned}\label{eq:relaxIE}
%&f^{(1)}_t + c f^{(1)} _x =\frac{1}{\epsilon} \left(E_1(\rho,m) - f^{(1)}\right),\\
%&f^{(2)}_t - c f^{(2)} _x =\frac{1}{\epsilon} \left(E_2(\rho,m) - f^{(2)}\right),\\
%&f^{(3)}_t  =\frac{1}{\epsilon} \left(E_3(\rho,m) - f^{(3)}\right),
%\end{aligned}
%\end{align}
%Where the equilibrium quantities are defined as follows
%\begin{align*}
%E_1(\rho,m,z) &=\frac{1}{2}F(\rho,m)+\frac{m}{2c},\\E_2(\rho,m,z) &=\frac{1}{2}F(\rho,m)-\frac{m}{2c},\\E_3(\rho,m,z) &=-F(\rho,m)+\rho
%\end{align*}
To this end we compute the adjoint equation system associated to \eqref{eq:relaxIE}, and equivalently to \eqref{BGK adjoint} we obtain the following
 \begin{equation}
\begin{aligned}\label{eq:relaxIEadjoint}
- & \lambda^{(1)}_t - c\lambda^{(1)}_x = - \frac{1}\epsilon \left( \lambda^{(1)} - \sum_k \lambda^{(k)}\left(  \partial_\rho E_k(\rho,m)+c \partial_m E_k(\rho,m)\right) \right), \\
- & \lambda^{(2)}_t + c\lambda^{(2)}_x = - \frac{1}\epsilon \left( \lambda^{(2)} - \sum_k \lambda^{(k)}\left(  \partial_\rho E_k(\rho,m)-c \partial_m E_k(\rho,m)\right) \right), \\
- & \lambda^{(3)}_t  = - \frac{1}\epsilon \left( \lambda^{(3)} - \sum_k \lambda^{(k)}\partial_\rho E_k(\rho,m) \right), 
 \end{aligned}
 \end{equation}
 complemented by the terminal conditions
 \[
 \lambda^{(1)}(T,x)  =  \partial_\rho J(\rho,m)+c\partial_m J(\rho,m),\,
 \lambda^{(2)}(T,x)  =  \partial_\rho J(\rho,m)-c\partial_m J(\rho,m),\, \lambda^{(3)}(T,x)  =  \partial_\rho J(\rho,m).
 \]
 
 We set up the control problem \eqref{eq:IE}--\eqref{eq:JIE} defining as reference density, and momentum the final state of system \eqref{eq:IE}  at time $T_f=0.15$ provided the following initial data
 \begin{equation}
 \begin{aligned}\label{realID}
 &\rho_d(0,x) = 1, \quad x\in[-2.5,2.5],
 \quad
 &m_d(0,x) = \begin{cases}\sin(\pi x), \quad  x\in [-1,1]\\ 0\qquad\qquad \textrm{otherwise}.\end{cases}
 \end{aligned}
 \end{equation}
 and zero flux boundary conditions.
 
In order to solve numerically problem \eqref{eq:relaxIE} --\eqref{eq:relaxIEadjoint}, we fix the relaxation parameter $\epsilon=0.01$. We discretize the space domain with an uniform grid of $N_x=320$ points, and with time step $\Delta t = 0.01$. In order to reduce the total variation of the initial data $(\rho^{(k)}_0,m^{(k)}_0)$ we introduce a filter $\mathcal{F}(\cdot)$ following the strategy proposed in \cite{HKK15}.  The optimization step is initialized using as starting guess the following data
\begin{equation}
 \begin{aligned}
 \rho^{(0)}(0,x) = 1,\qquad  m^{(0)}(0,x) = 0, \quad x\in[-2.5,2.5].
 \end{aligned}
 \end{equation} 
Then at each iteration $k=0,1,\ldots$ the initial data $\rho_0^{(k)}(x),m_0^{(k)}(x)$ is updated with gradient method with Barzilai-Borwein descent step, \cite{BB}.
%\begin{aligned}
%u^{(k+1)} _0=u_0^{(k)}-\sigma_k p^{(k)}, \qquad m\geq 0,
%\]

\begin{figure}[tb]
\centering
\includegraphics[scale = 0.35]{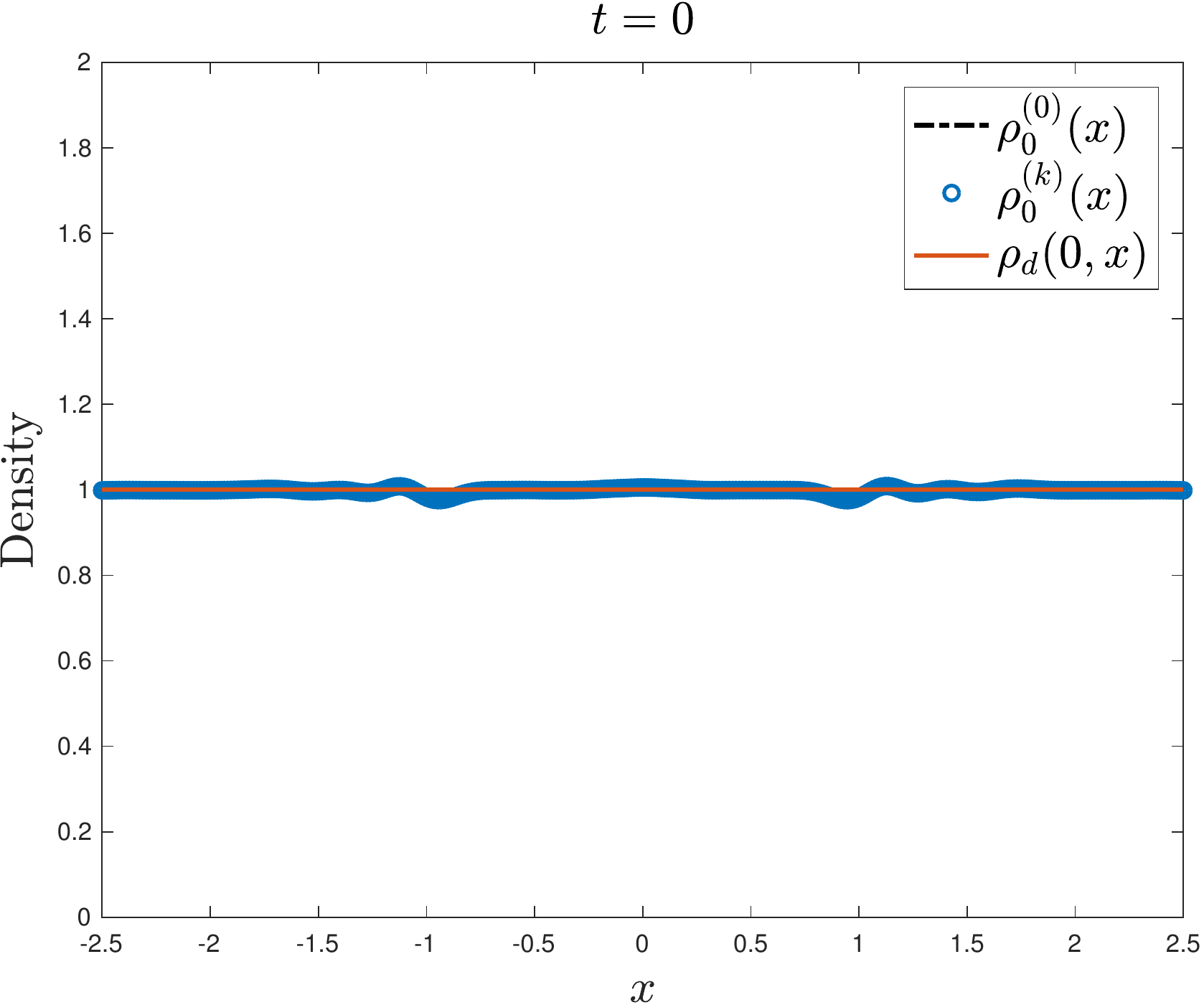}
\quad\qquad
\includegraphics[scale = 0.35]{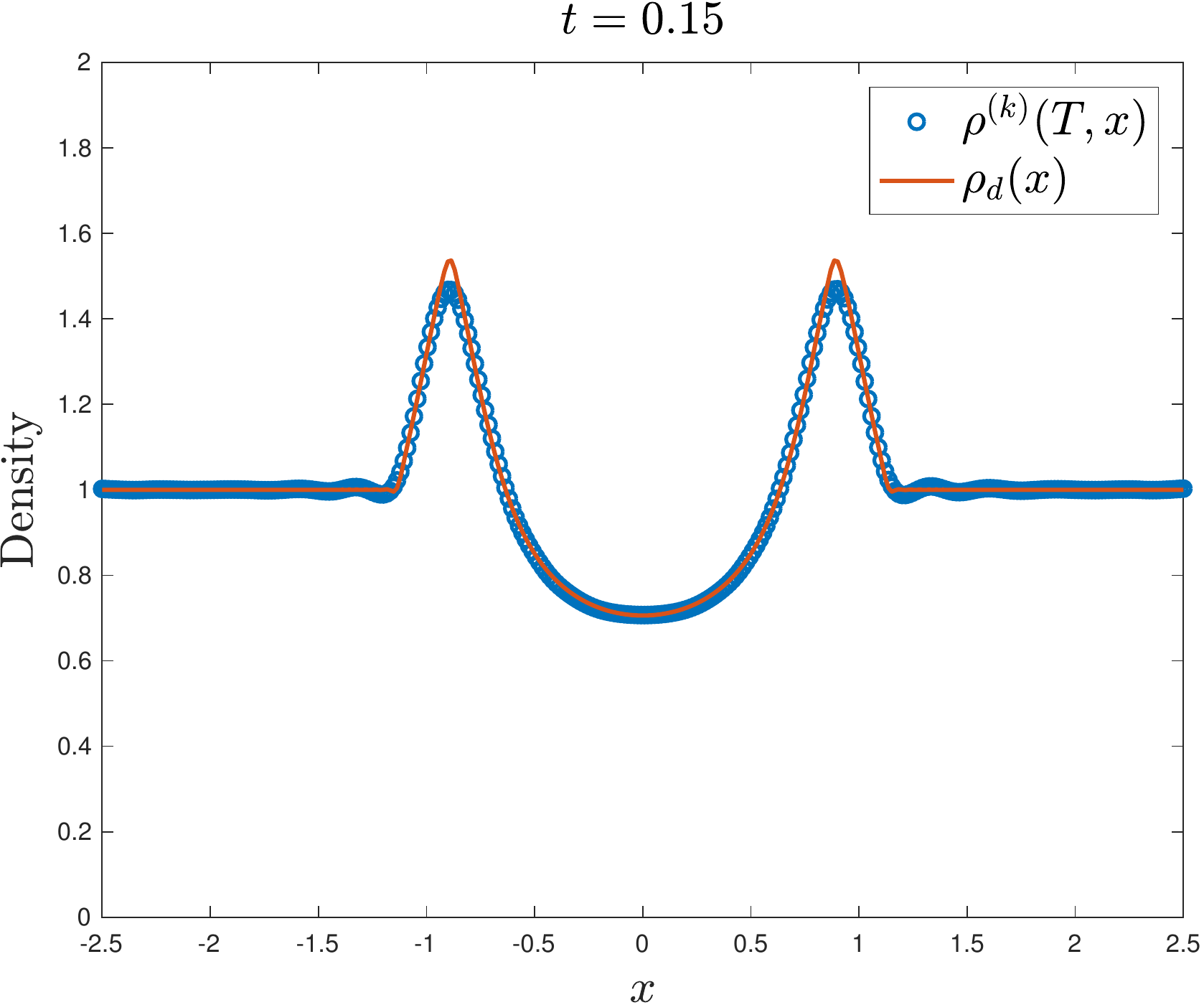}
\\
\includegraphics[scale = 0.35]{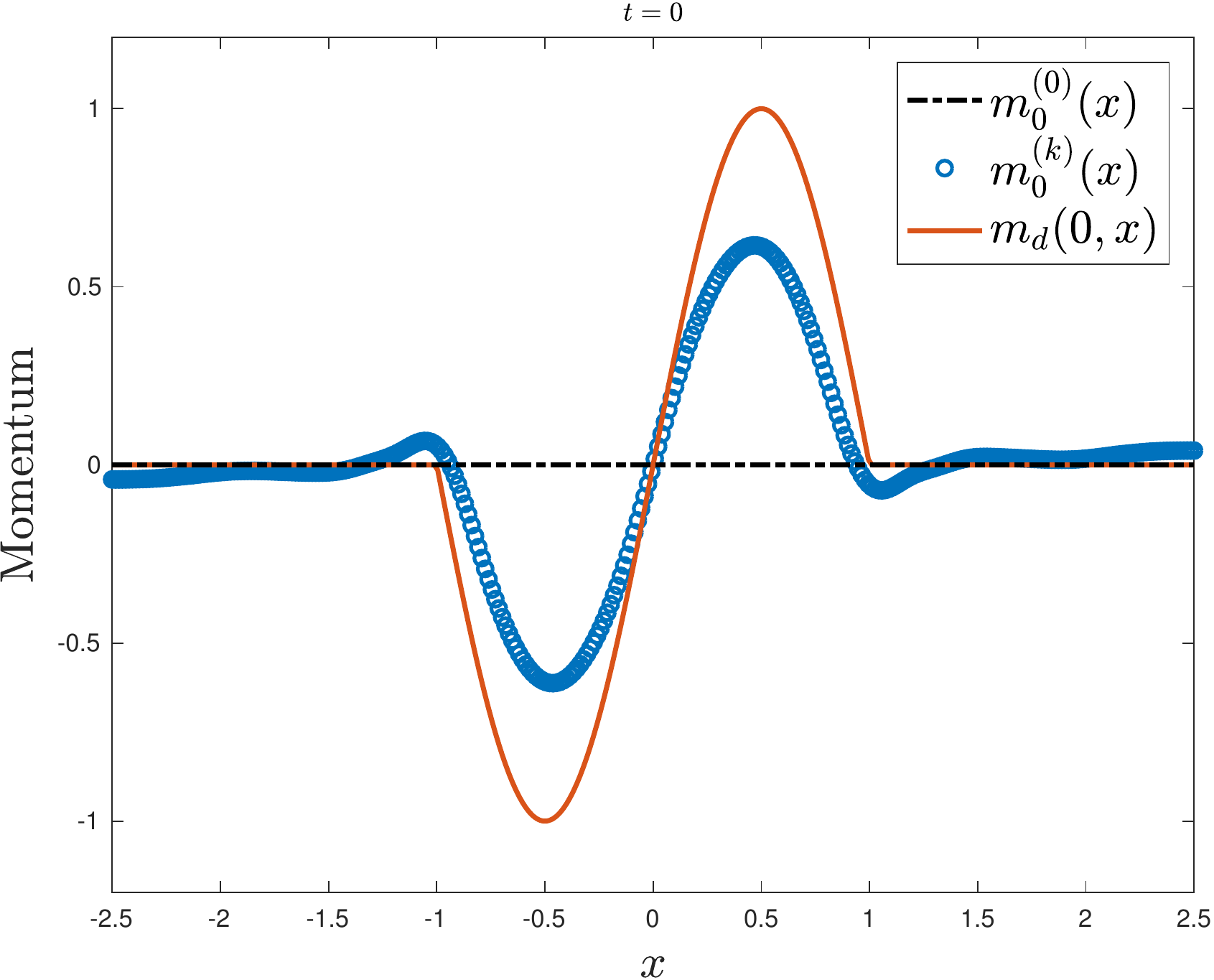}
\quad\qquad
\includegraphics[scale = 0.35]{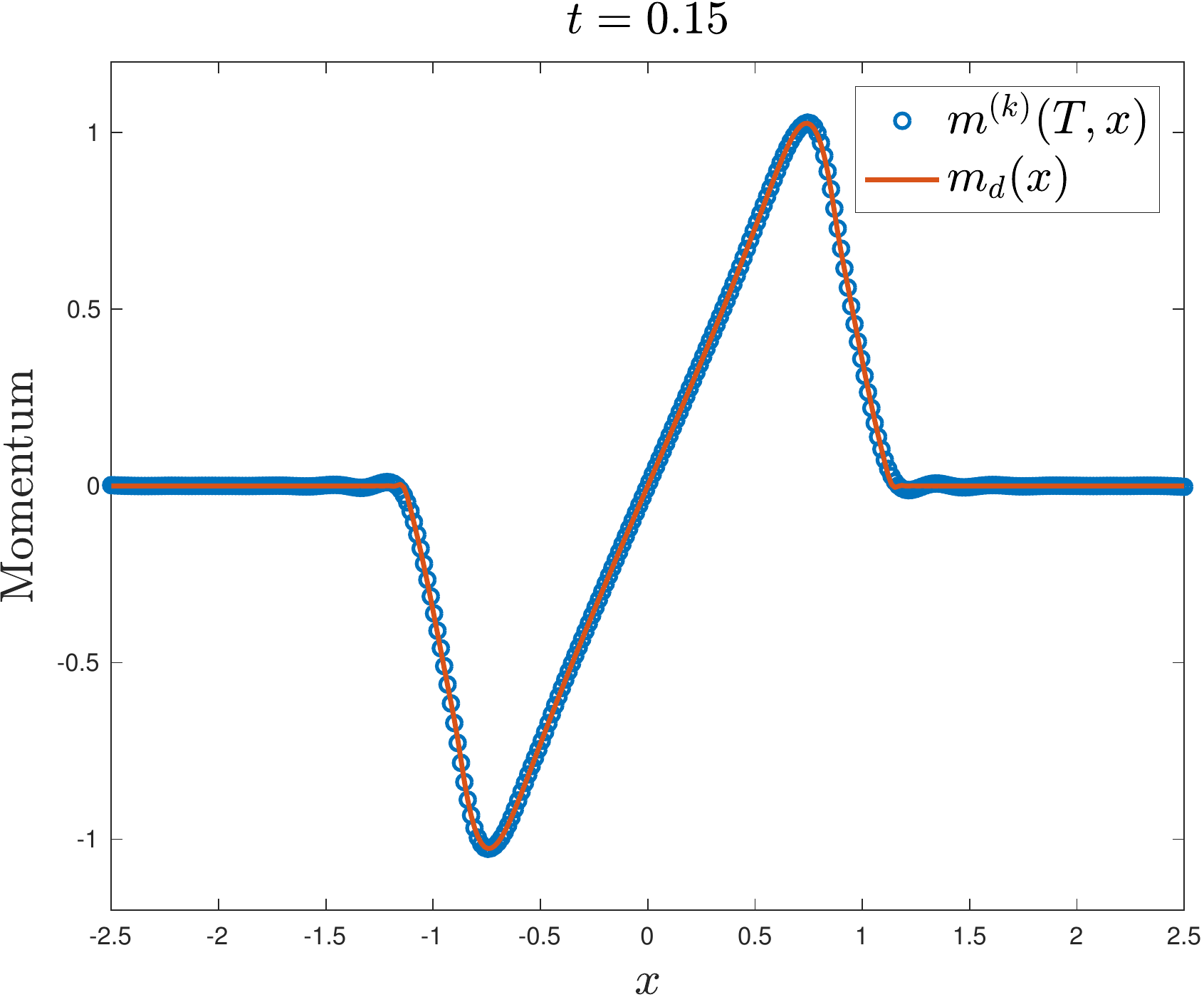}
\caption{{\em Broadwell model}. We consider $\epsilon=0.01$, using $N_x  =320$ space points and $\dt = 1\times 10^{-2}$. Top row represents the initial and final time of the density $\rho(t,x)$, comparing the optimal control $\rho^{(k)}$ with respect to the reference $\rho_d$ at initial (left plot) and final time (right plot). Bottom row compares the momentum $m^{(k)}(t,x)$ at initial (left plot) and final time (right plot)  with respect to the reference solution $m_d(x)$.
}\label{Fig005}
\end{figure}

\begin{figure}[tb]
\centering
\includegraphics[scale = 0.45]{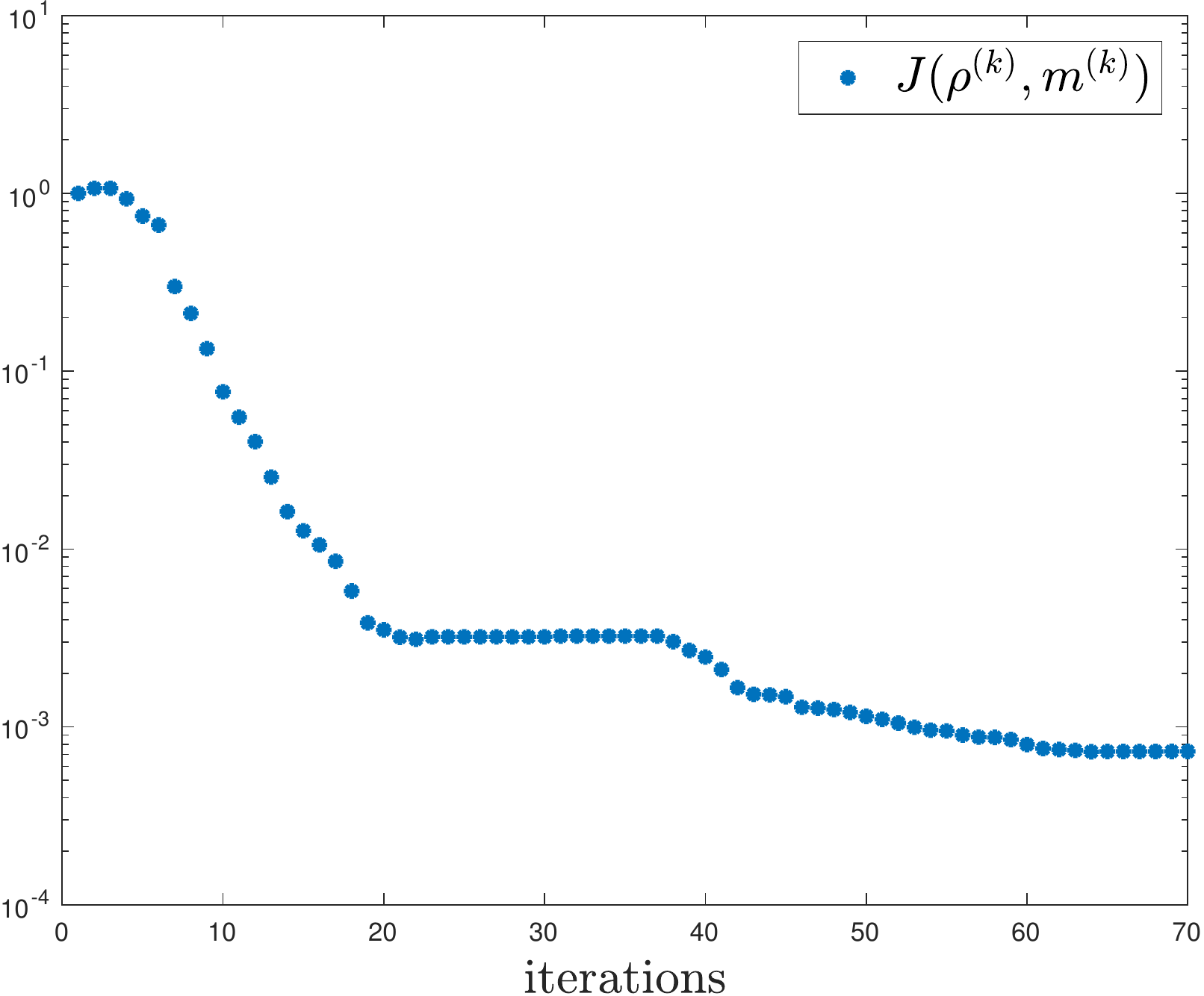}
\caption{{\em Broadwell model.} Decrease of the functional $J(\rho,m)$ evaluated in $\rho^{(k)},m^{(k)}$, at each iteration $k$ of the optimization process.}\label{Fig006}
\end{figure}

 We report in Figure \ref{Fig005} the result of the optimization process. Top row depicts the evolution of the density, whereas bottom row refers to momentum evolution. On the left-hand side the initial value $(\rho^{(0)}_0(x), m_0^{(0)}(x))$ is compared  with the control $(\rho_0^{(k)}(x),m_0^{(k)}(x))$ obtained after $k=70$ iterations of the optimization process, and the true initial data defined by \eqref{realID}. The right-hand side column depicts the density and momentum at final time $T= 0.15$ comparing the reference solution $(\rho_d(x),m_d(x))$ with respect to $(\rho^{(k)}(T,x),m^{(k)}(T,x))$.
 Finally Figure \ref{Fig006} reports the decrease of the functional $J(\rho,m)$ evaluated at each iteration of the optimization process.

%\begin{figure}[tb]
%\centering
%\includegraphics[scale = 0.35]{datarho0_broadwell.pdf}
%\quad\qquad
%\includegraphics[scale = 0.35]{datarhoTf_broadwell.pdf}
%\caption{The BDF(2) integration of the system of \eqref{implement adjoint} has been implemented for the linear transport $F(u)=u.$ Two velocities are considered,  with $N_x  =640$ space points and $\dt = 4.47127\times 10^{-3}$. From left to right we show different values of $\epsilon$, with $\epsilon\in{\{1,0.1,0.0001\}}$. In the top row the terminal data $p_T(x)$, $T=1$, as well as the numerical result at initial time $p(0,x)$ are reported for the different values of $\epsilon$. Bottom row depict the density $p(t,x) = \lambda^1(t,x)+\lambda^2(t,x)$ for the different value of $\epsilon$. Note that for small $\epsilon$ the pure transport equation is obtained.}\label{Fig002}
%\end{figure}
%

\section{Conclusion}

We analyze linear multi-step schemes for control problems of ordinary differential equations and hyperbolic balance laws. In the case of ordinary differential equations we show theoretically and numerically that only BDF methods are  consistent discretization of the corresponding optimality systems up to high--order. 
The BDF methods may also be used as higher order discretization of
relaxation systems in combination with a Lagrangian scheme. We derive the corresponding adjoint equations and we show that this system can again be discretized by a BDF type method. The numerically observed convergence rates confirm the expected behavior both for ordinary differential systems, as well as hyperbolic balance laws.

\appendix
\section{ Definition of BDF, Adams--Moulton and Adams--Bashfort Formulas}\label{tables}

In view of the scheme \eqref{scheme} each scheme is represented by 
two vectors $a,b$ with $a=(a_0,\dots,a_{s-1}) \in \R^s$ and $b=(b_{-1},b_0,\dots, b_{s-1}) 
\in \R^{s+1}$ for an $s-$stage
scheme. Only in the case of BDF schemes we have $b \in \R^{s+1},$ otherwise
we have $b\in \R^s.$ For the schemes implemented in this paper we use the 
following schemes. 

\begin{table}[htb]
\center
\begin{tabular}{llll}
Name 		  &  s & $a^t$ & $b^t$ \\ \hline
Implicit Euler  (BDF(1)) & 1 & -1& (1,0) \\
Explicit Euler & 1 & 0 & (0,1) \\
\hline
{\em BDF methods} \\ 
BDF(2) & 2 & (-4/3,1/3) & (2/3,0,0) \\
BDF(3) & 3 & (-18/11,9/11,-2/11) & (6/11,0,0,0) \\
BDF(4) & 4 & (-48/25,36/25,-16/25,3/25) & ( 12/25,0,0,0,0) \\
\hline
{\em Adams--Bashfort (AB) methods } \\
AB(2) & 2 & (-1,0) & (0,3/2,-1/2) \\
AB(3) & 3 & (-1,0,0) & (0,23/12,-4/3,5/12) \\
\hline
{\em Adams--Moulton (AM) methods } \\
AM(4) & 4 & (-1,0,0,0) & (251,646,-264,106,-19)/270\\
\hline
\end{tabular}
\end{table}

\subsubsection*{Acknowledgments}
This work has been supported by DFG HE5386/13,14,15-1, by 
 the DAAD--MIUR project, KI-Net  and by the INdAM-GNCS 2018 project {\it Numerical methods for multi-scale control problems and applications}.
 
%%%%%%%%%%%%%%%%%%%%%%%%%%%%%%%%%%%%%%%%%%%

\end{document}